\RequirePackage{etoolbox}
\csdef{input@path}{{style/}{graphics/}}
\documentclass[aos,MSNbibl,rotating,nameyear,seceqn,dvips]{arximspdf}
\usepackage{mathrsfs,mathbh,dcolumn}

\doi{10.1214/14-AOS1299}
\volume{43}
\issue{2}
\pubyear{2015}
\firstpage{741}
\lastpage{768}
\docsubty{FLA}

\makeatletter
\def\mid{|}
\newcolumntype{d}[1]{D{.}{.}{#1}}
\def\boldsymbol{\mathbf}
\newcommand{\rrVert}{\Vert}
\newcommand{\llVert}{\Vert}
\def\cal{\mathcal}

\def\II{II}
\def\III{III}
\def\IV{IV}
\def\B{\mathbb{B}}

\newproclaim{definition}{Definition}[section]

\newtheorem{lemma}{Lemma}[section]
\newtheorem{proposition}{Proposition}[section]
\newtheorem{theorem}{Theorem}[section]

\newproclaim{remark}{Remark}[section]
\newproclaim{example}{Example}[section]

\newcommand{\indicator}[1]{\mathbh{1}_{\{ {#1} \} }}
\makeatother

\begin{document}
\begin{frontmatter}

\title{Time-varying nonlinear regression models: Nonparametric
estimation and model selection}
\runtitle{Time-varying nonlinear regression models}

\begin{aug}
%
\begin{aug}
\author[A]{\fnms{Ting}~\snm{Zhang}\corref{}\thanksref{C1}\ead[label=e1]{tingz@bu.edu}}
\and
\author[B]{\fnms{Wei~Biao}~\snm{Wu}\thanksref{C2}\ead[label=e2]{wbwu@galton.uchicago.edu}}
\runauthor{T. Zhang and W.~B. Wu}
\thankstext{C1}{Supported in part by NSF Grant DMS-14-61796.}
\thankstext{C2}{Supported in part by NSF Grants DMS-14-05410 and DMS-11-06790.}

\affiliation{Boston University and University of Chicago}

\address[A]{Department of Mathematics and Statistics\\
Boston University\\
Boston, Massachusetts 02215\\
USA\\
\printead{e1}}

\address[B]{Department of Statistics\\
University of Chicago\\
Chicago, Illinois 60637\\
USA\\
\printead{e2}}
\end{aug}
\end{aug}

%
\received{\smonth{5} \syear{2014}}
%
\revised{\smonth{12} \syear{2014}}

%
\begin{abstract}
This paper considers a general class of nonparametric time series
regression models where the regression function can be time-dependent.
We establish an asymptotic theory for estimates of the time-varying
regression functions. For this general class of models, an important
issue in practice is to address the necessity of modeling the
regression function as nonlinear and time-varying. To tackle this, we
propose an information criterion and prove its selection consistency
property. The results are applied to the U.S. Treasury interest rate data.
\end{abstract}

%
\begin{keyword}[class=AMS]
\kwd[Primary ]{62G05}
\kwd{62G08}
\kwd[; secondary ]{62G20}
\end{keyword}

\begin{keyword}
\kwd{Information criterion}
\kwd{local linear estimation}
\kwd{nonparametric model selection}
\kwd{nonstationary processes}
\kwd{time-varying nonlinear regression models}
\end{keyword}

\end{frontmatter}
%
\section{Introduction}\label{secintro}
Consider the time-varying regression model
%
%
\begin{equation}
\label{eqnformmxt} \mbox{Model I:}\quad y_i = m_i(\boldsymbol
x_i) + e_i,\qquad i = 1,\ldots,n,
\end{equation}
where $y_i$, $\boldsymbol x_i$ and $e_i$ are the responses, the
predictors and the errors, respectively, and $m_i(\cdot) = m(\cdot,
i/n)$ is a time-varying regression function. Here $m\dvtx\mathbb R^d
\times[0,1] \to\mathbb R$ is a smooth function, and $i/n$, $i =
1,\ldots,n$, represents the time rescaled to the unit interval. Model
$\mathrm I$ is very general. If $m_i(\cdot)$ is not time-varying, then
(\ref{eqnformmxt}) becomes
\[
\mbox{Model \II:}\quad y_i = \mu(\boldsymbol x_i) +
e_i,\qquad i = 1,\ldots,n.
\]
Model $\mathrm{\II}$ has been extensively studied in the literature;
see \citet{Robinson1983}, \citet{GyorfiHardleSardaVieu1989}, \citet
{FanYao2003} and \citet{LiRacine2007}, among others. As an
important example, (\ref{eqnformmxt}) can be viewed as the discretized
version of the nonstationary diffusion process
%
%
\begin{equation}
\label{eqnSDExt} dy_t = m(y_t,t/T)\,dt +
\sigma(y_t,t/T)\,d\B_t,
\end{equation}
where $\{\B_s\}_{s \in\mathbb{R}}$ is a standard Brownian motion,
$m(\cdot,\cdot)$ and $\sigma(\cdot,\cdot)$ are, respectively, the drift
and the volatility functions, which can both be time-varying, and $T$
is the time horizon under consideration. If the functions $m(\cdot
,\cdot
)$ and $\sigma(\cdot,\cdot)$ do not depend on time, then (\ref
{eqnSDExt}) becomes the stationary diffusion process
%
%
\begin{equation}
\label{eqnSDEx} dy_t = \mu(y_t)\,dt + \gamma(y_t)\,d
\B_t,
\end{equation}
which relates to model $\mathrm{\II}$. There is a huge literature on
modeling interest rates data by (\ref{eqnSDEx}). For example, \citet
{Vasicek1977} considered model (\ref{eqnSDEx}) with linear drift
function $\mu(x) = \beta_0 + \beta_1 x$ and constant volatility
$\gamma
(x) \equiv\gamma$, where $\beta_0,\beta_1,\gamma$ are unknown
parameters. \citet{Courtadon1982}, \citet{CoxIngersollRoss1985} and
\citet{ChanKarolyiLongstaffSanders1992} considered nonconstant
volatility functions. \citet{AitSahalia1996}, \citet{Stanton1997} and
\citet{LiuWu2010} studied model (\ref{eqnSDEx}) with nonlinear drift
functions. See \citet{Zhao2008} for a review. However, due to policy
and societal changes, those models with static relationship between
responses and predictors may not be suitable. Here we shall study
estimates of time-varying regression function $m_i(\cdot)$ for model
(\ref{eqnformmxt}).

For model $\mathrm{\II}$, let $K_S(\cdot)$ be a $d$-dimensional
kernel function
%
%
\begin{equation}\quad
\label{eqnftilde} \tilde T_n(\boldsymbol u) = \frac{1 }{ nh_n^d} \sum
_{i=1}^n y_i K_S
\biggl(\frac{\boldsymbol u-\boldsymbol x_i }{ h_n} \biggr),\qquad \tilde
f_n(\boldsymbol u) =
\frac{1 }{ nh_n^d} \sum_{i=1}^n
K_S \biggl(\frac{\boldsymbol u-\boldsymbol x_i }{ h_n} \biggr),
\end{equation}
where $h_n$ be a bandwidth sequence. We can then apply the traditional
Nadaraya--Watson estimate for the regression function $\mu(\cdot)$,
%
%
\begin{equation}
\label{eqnmuhat} \hat\mu_n(\boldsymbol u) = \frac{\tilde
T_n(\boldsymbol u) }{\tilde
f_n(\boldsymbol u)},\qquad
\boldsymbol u \in\mathbb R^{d}.
\end{equation}
If the process $(\boldsymbol x_i)$ is stationary, then $\tilde f_n$ is
the kernel density estimate of its marginal density. For stationary
processes, an asymptotic theory for these nonparametric estimators has
been developed by many researchers, including \citet{Robinson1983},
\citet{CastellanaLeadbetter1986}, \citet{Silverman1986}, \citet
{GyorfiHardleSardaVieu1989}, \citet{Yu1993}, \citet
{Tjotheim1994}, \citet{WandJones1995}, \citet{Bosq1996}, \citet
{Neumann1998}, \citet{NeumannKreiss1998}, \citet{FanYao2003} and
\citet{LiRacine2007}, among others. However, the case of
nonstationary processes has been rarely touched. \citet
{HallMullerWu2006} considered the situation that the underlying
distribution evolves with time and proposed a nonparametric
time-dynamic density estimator. Assuming independence, they proved the
consistency of their kernel-type estimators and applied the results to
fast mode tracking. Following the spirit of \citet
{HallMullerWu2006}, \citet{Vogt2012} considered a kernel estimator
of the time-varying regression model (\ref{eqnformmxt}), and
established its asymptotic normality and uniform bound under the
classical strong mixing conditions. In Sections~\ref{subsecCLTs} and
\ref{subsecunifbnds}, we advance the nonparametric estimation theory
for the time-varying regression model (\ref{eqnformmxt}) under the
framework of \citet{DraghicescuGuillasWu2009}, which is convenient
to use and often leads to optimal asymptotic results.

Apart from model $\mathrm{\II}$, model $\mathrm I$ contains another
important special case: the time-varying coefficient linear regression model
\[
\mbox{Model \III:} \quad y_i = \boldsymbol x_i^\top
\bolds\beta_i + e_i,\qquad i = 1,\ldots,n,
\]
where $^\top$ is the transpose and $\bolds\beta_i = \bolds\beta(i/n)$ for some smooth function $\bolds\beta\dvtx [0,1]
\to
\mathbb R^d$. The traditional linear regression model
\[
\mbox{Model \IV:} \quad y_i = \boldsymbol x_i^\top
\bolds\theta+ e_i,\qquad  i = 1,\ldots,n,
\]
where $\bolds\theta\in\mathbb R^d$ is the regression
coefficient, is a special case of model $\mathrm{\III}$. Estimation of
$\bolds\beta(\cdot)$ has been considered by \citet
{HooverRiceWuYang1998}, \citeauthor{FanZhangJT2000}
(\citeyear{FanZhangJT2000,FanZhangWY2000}),
\citet{HuangWuZhou2004}, \citet
{RamsaySilverman2005}, \citet{Cai2007} and \citet{ZhouWu2010},
among others. The problem of distinguishing between models $\mathrm
{\III
}$ and $\mathrm{\IV}$ has been studied in the literature mainly by
means of hypothesis testings; see, for example, \citet{Chow1960},
\citet{BrownDurbinEvans1975}, \citet{NabeyaTanaka1988}, \citet
{LeybourneMcCabe1989}, \citet{Nyblom1989}, \citet
{PlobergerKramerKontrus1989}, \citet{Andrews1993}, \citet
{DavisHuangYao1995}, \citet{LinTerasvirta1999} and \citet
{HeTerasvirtaGonzalez2009}. On the other hand, model $\mathrm{\IV}$
specifies a linear relationship upon model $\mathrm{\II}$, and there is
a huge literature on testing parametric forms of $\mu(\cdot)$; see
\citet{AzzaliniBowman1993}, \citet{GonzalezCao1993}, \citet
{HardleMammen1993}, \citet{Zheng1996}, \citet{Dette1999}, \citet
{FanZhangZhang2001}, \citet{ZhangDette2004} and \citet
{ZhangWu2011}, among others. Nevertheless, model selection between
models $\mathrm{\II}$ and $\mathrm{\III}$ received much less attention.
Note that both of them are nested in the general model $\mathrm I$, and
they all cover the linear regression model $\mathrm{\IV}$. It is
desirable to develop a model selection criterion. An information
criterion is proposed in Section~\ref{subsecModelSelection}, where its
consistency property is obtained.

The rest of the paper is organized as follows. Section~\ref
{secframework} introduces the model setting. Main results are stated
in Section~\ref{secmain} and are proved in Section~\ref{secappendix}
with some of the proofs postponed to the supplementary material [\citet
{ZhangWu2014}]. A simulation study is given in Section~\ref
{secnumerexp} along with an application to the U.S. Treasury interest
rate data.

\section{Model setting}
\label{secframework}
For estimation of model $\mathrm I$, temporal dynamics should be taken
into consideration. Let $K_T(\cdot)$ be a temporal kernel function
(kernel function for time), $b_n$ be another sequence of bandwidths and
$w_{b_n,i}(t) = K_T \{(i/n-t) / b_n\} \{S_2(t)-(t-i/n)S_1(t)\} / \{
S_2(t)S_0(t)-S_1^2(t)\}$ be the local linear weights, where $S_l(t) =
\sum_{j=1}^n (t-j/n)^l K_T\{(j/n-t)/b_n\}$, $l \in\{0,1,2\}$. Let
$K_{S,h_n}(\cdot) = h_n^{-d}K_S(\cdot/h_n)$,
%
%
\begin{eqnarray}
\label{eqnfhat} \hat f_n(\boldsymbol u,t)& =& \sum
_{i=1}^n K_{S,h_n}(\boldsymbol u-\boldsymbol
x_i)w_{b_n,i}(t),
\nonumber
\\[-8pt]
\\[-8pt]
\nonumber
\hat T_n(\boldsymbol u,t) &=&
\sum_{i=1}^n y_i
K_{S,h_n}(\boldsymbol u-\boldsymbol x_i)w_{b_n,i}(t),
\end{eqnarray}
we consider the time-varying kernel regression estimator
%
%
\begin{equation}
\label{eqnmhat} \hat m_n(\boldsymbol u,t) = \frac{\hat
T_n(\boldsymbol u,t) }{\hat
f_n(\boldsymbol u,t)}.
\end{equation}
\citet{HallMullerWu2006} proved the uniform consistency of $\hat
f_n$ in (\ref{eqnfhat}) by assuming that $\boldsymbol x_1, \ldots,
\boldsymbol x_n$ are independent. To allow nonstationary and dependent
observations, we assume
%
%
\begin{equation}
\label{eqnframeworkx} \boldsymbol x_i = \boldsymbol G(i/n;
\bolds{\mathcal{ H}}_i), \qquad\mbox{where } \bolds{\mathcal{ H}}_i
= (\ldots, \bolds\xi_{i-1}, \bolds\xi_i)
\end{equation}
and $\bolds\xi_k$, $k \in\mathbb Z$, are independent and
identically distributed (i.i.d.) random vectors, and $\boldsymbol G$ is
a measurable function such that $\boldsymbol G(t;
\bolds{\mathcal{ H}}_i)$ is well defined for each $t \in[0,1]$. Following \citet
{DraghicescuGuillasWu2009}, the framework (\ref{eqnframeworkx})
suggests locally strict stationarity and is convenient for asymptotic
study. For the error process, we assume that
%
%
\begin{equation}
\label{eqndefeiiideta} e_i = \sigma_i(\boldsymbol
x_i)\eta_i = \sigma(\boldsymbol x_i, i/n)
\eta_i,
\end{equation}
where $\sigma(\cdot, \cdot)\dvtx\mathbb R^d \times[0,1] \to
\mathbb R$ is
a smooth function, and $(\eta_i)$ is a sequence of random variables
satisfying $E(\eta_i \mid\boldsymbol x_i) = 0$ and $E(\eta_i^2 \mid
\boldsymbol x_i) = 1$. At the outset (cf. Sections~\ref
{subsecCLTs}--\ref{subsecModelSelection}) we assume that $\eta_k$, $k
\in\mathbb Z$, are i.i.d. and independent of
$\bolds{\mathcal{H}}_j$, $j \in\mathbb Z$. The latter assumption can be relaxed (though
technically much more tedious) to allow models with correlated errors
and nonlinear autoregressive processes; see Section~\ref{subsecextensions}.

For a random vector $\boldsymbol Z$, we write $\boldsymbol Z \in
\mathcal L^q$, $q > 0$ if $\|\boldsymbol Z\| = \{E(|\boldsymbol Z|^q)\}
^{1/q} < \infty$ where $|\cdot|$ is the Euclidean vector norm, and we
denote $\| \cdot\| = \| \cdot\|_2$. Let $F_1(\boldsymbol u,t \mid
\bolds{\mathcal{H}}_k) = \mathrm{pr}\{\boldsymbol G(t;\bolds{\mathcal{ H}}_{k+1}) \leq\boldsymbol u
\mid\bolds{\mathcal{ H}}_k\}
$ be the one-step ahead predictive or conditional distribution function
and $f_1(\boldsymbol u,t \mid\bolds{\mathcal{ H}}_k) = \partial^d
F_1(\boldsymbol u,t \mid\bolds{\mathcal{H}}_k) / \partial
\boldsymbol u$ be the corresponding conditional density. Let
$(\bolds\xi_i')$ be an i.i.d. copy of $(\bolds\xi_j)$ and
$\bolds{\mathcal{ H}}_k' = (\ldots,\bolds\xi
_{-1},\bolds\xi_0',\bolds\xi_1,\ldots,\bolds\xi_k)$ be the coupled
shift process. We define the predictive dependence measure
%
%
\begin{equation}
\label{eqnpsikq} \psi_{k,q} = \sup_{t \in[0,1]} \sup
_{\boldsymbol u \in\mathbb R^d} \bigl\| f_1(\boldsymbol u,t \mid
\bolds{\mathcal{H}}_k) - f_1\bigl(\boldsymbol u,t \mid\bolds{\mathcal{ H}}_k'\bigr)\bigr\|_q.
\end{equation}
Quantity (\ref{eqnpsikq}) measures the contribution of $\bolds\xi
_0$, the innovation at step 0, on the conditional or predictive
distribution at step $k$. We shall make the following assumptions:
\begin{longlist}[(A1)]
\item[(A1)] smoothness (third order continuous differentiability):
$f,m,\break \sigma\in\mathcal C^3(\mathbb R^d \times[0, 1])$;
\item[(A2)] short-range dependence: $\Psi_{0,2} < \infty$, where
$\Psi
_{m,q} = \sum_{k = m}^\infty\psi_{k,q}$;
\item[(A3)] there exists a constant $c_0 < \infty$ such that almost surely,
\[
\sup_{t \in[0,1]} \sup_{\boldsymbol u \in\mathbb{R}^d} \bigl\{
f_1(\boldsymbol u,t \mid\bolds{\mathcal{ H}}_0) + \bigl|
\partial^d f_1(\boldsymbol u,t \mid\bolds{\mathcal{H}}_0)/\partial\boldsymbol u\bigr|\bigr\} \leq c_0.
\]
\end{longlist}
Condition (A3) implies that the marginal density $f(\boldsymbol u,t) =
E\{f_1(\boldsymbol u,t \mid\bolds{\mathcal{ H}}_0)\} \le c_0$.\vadjust{\goodbreak}

\section{Main results}\label{secmain}
\subsection{Nonparametric kernel estimation}\label{subsecCLTs}
Throughout the paper, we assume that the kernel functions $K_S(\cdot)$
and $K_T(\cdot)$ are both symmetric and twice continuously
differentiable on their support $[-1,1]^d$ and $[-1,1]$, respectively,
and $\int_{[-1,1]^d} K_S(\boldsymbol s) \,d\boldsymbol s = \int_{-1}^1
K_T(v) \,dv = 1$.\vspace*{1pt} Denote by ``$\Rightarrow$'' convergence in
distribution. Theorem~\ref{thmCLTfhatmhat} provides the asymptotic
normality of the time-varying kernel estimators (\ref{eqnfhat}) and
(\ref{eqnmhat}), while Theorem~\ref{thmCLTftildemuhat} concerns the
time-constant estimators (\ref{eqnftilde}) and (\ref{eqnmuhat}).

%
\begin{theorem}\label{thmCLTfhatmhat}
Assume \textup{(A1)--(A3)} and $\eta_i \in\mathcal L^p$, $p > 2$ are i.i.d. Let
$(\boldsymbol u,t) \in\mathbb{R}^d \times(0,1)$ be a fixed point. If
$b_n \to0$, $h_n \to0$ and $nb_nh_n^d \to\infty$, then
%
%
\begin{equation}
\label{eqnCLTfhat} \bigl(nb_nh_n^d
\bigr)^{1/2} \bigl[\hat f_n(\boldsymbol u,t) - E\bigl\{\hat
f_n(\boldsymbol u,t)\bigr\}\bigr] \Rightarrow N \bigl\{0, f(
\boldsymbol u,t) \lambda_{K_S} \lambda_{K_T}\bigr\},
\end{equation}
where $\lambda_{K_T} = \int_{-1}^1 K_T(v)^2 \,dv$ and $\lambda_{K_S} =
\int_{[-1,1]^d} K_S(\boldsymbol s)^2 \,d\boldsymbol s$. If in addition\break
$f(\boldsymbol u,t) > 0$, then
%
%
\begin{equation}
\label{eqnCLTmhat} \bigl(nb_nh_n^d
\bigr)^{1/2} \biggl[\hat m_n(\boldsymbol u,t) -
\frac{E\{\hat
T_n(\boldsymbol u,t)\} }{ E\{\hat f_n(\boldsymbol u,t)\}} \biggr]
\Rightarrow N \biggl\{0, \frac{\sigma(\boldsymbol u,t)^2 \lambda
_{K_S} \lambda
_{K_T} }{ f(\boldsymbol u,t)} \biggr\}.
\end{equation}
\end{theorem}

Let $H_f(\boldsymbol u,t) = \{\partial^2 f(\boldsymbol u,t) / \partial
u_i\, \partial u_j\}_{1 \leq i,j \leq d}$\vspace*{1pt} be the Hessian matrix of the
density function $f$ with respect to $\boldsymbol u$. Denote
$f^{(\boldsymbol 0,2)}(\boldsymbol u,t) = \partial^2 f(\boldsymbol u,t)
/ \partial t^2$, and we use the same notation for the product function
$(mf)(\boldsymbol u,t) = m(\boldsymbol u,t) f(\boldsymbol u,t)$. Then
for any point $(\boldsymbol u,t) \in\mathbb{R}^d \times(0,1)$ with
$f(\boldsymbol u,t) > 0$, we have
\[
E\bigl\{\hat f_n(\boldsymbol u,t)\bigr\} = f(\boldsymbol u,t) +
\frac{h_n^2 }{2} \mathrm{tr} \bigl\{H_f(\boldsymbol u,t)\bolds\kappa_S \bigr\} + \frac{b_n^2 }{2}
 f^{(\boldsymbol
0,2)}(\boldsymbol
u,t) \kappa_T + O\bigl(b_n^3 +
h_n^3\bigr),
\]
where $\mathrm{tr}(\cdot)$ is the trace operator
\[
\bolds\kappa_S = \int_{[-1,1]^d}
K_S(\boldsymbol s) \boldsymbol s \boldsymbol s^\top\, d
\boldsymbol s,\qquad \kappa_T = \int_{-1}^1
K_T(v) v^2 \,dv
\]
and
\begin{eqnarray*}
\frac{E\{\hat T_n(\boldsymbol u,t)\} }{ E\{\hat f_n(\boldsymbol u,t)\}
} & = & m(\boldsymbol u,t) + \frac{h_n^2 }{2 f(\boldsymbol u,t)}
\mathrm{tr}
\bigl[\bigl\{H_{mf}(\boldsymbol u,t) - m(\boldsymbol
u,t)H_f(\boldsymbol u,t)\bigr\} \bolds\kappa_S
\bigr]
\\
& &{} + \frac{b_n^2 }{2f(\boldsymbol u,t)} \bigl\{(mf)^{(\boldsymbol
0,2)}(\boldsymbol u,t) - m(
\boldsymbol u,t) f^{(\boldsymbol
0,2)}(\boldsymbol u,t)\bigr\} \kappa_T
\\
& &{} + O\bigl(b_n^3 + h_n^3
\bigr).
\end{eqnarray*}
Hence (\ref{eqnfhat}) and (\ref{eqnmhat}) are consistent estimates of
the local density function $f$ and the regression function $m$,
respectively. The asymptotic mean squared error (AMSE) optimal
bandwidths satisfy $b_n \asymp n^{-1/(d+5)}$ and $h_n \asymp
n^{-1/(d+5)}$. Here for positive sequences $(s_n)$ and $(r_n)$, we
write $s_n \asymp r_n$ if $s_n/r_n + r_n/s_n$ is bounded for all large $n$.

%
\begin{theorem}\label{thmCLTftildemuhat}
Assume \textup{(A1)--(A3)} and $\eta_i \in\mathcal L^p$, $p > 2$. If $h_n \to
0$ and $nh_n^d \to\infty$, then
%
%
\begin{equation}
\label{eqnCLTftilde} \bigl(nh_n^d\bigr)^{1/2} \bigl[
\tilde f_n(\boldsymbol u) - E\bigl\{\tilde f_n(\boldsymbol
u)\bigr\}\bigr] \Rightarrow N \bigl\{0, \bar f(\boldsymbol u)
\lambda_{K_S}\bigr\}, \qquad \boldsymbol u \in\mathbb{R}^d,
\end{equation}
where $\bar f(\boldsymbol u) = \int_0^1 f(\boldsymbol u,t) \,d t$. If in
addition $\bar f(\boldsymbol u) > 0$, then
%
%
\begin{equation}
\label{eqnCLTmuhat} \bigl(nh_n^d\bigr)^{1/2} \biggl[
\hat\mu_n(\boldsymbol u) - \frac{E\{\tilde
T_n(\boldsymbol u)\} }{ E\{\tilde f_n(\boldsymbol u)\}} \biggr]
\Rightarrow N
\bigl\{0, \tilde V(\boldsymbol u) \lambda_{K_S}\bigr\},
\end{equation}
where, letting $\bar m(\boldsymbol u) = \int_0^1 m(\boldsymbol u,t)
f(\boldsymbol u,t) \,dt / \bar f(\boldsymbol u)$, the variance function
\[
\tilde V(\boldsymbol u) = \bar f(\boldsymbol u)^{-2} \int
_0^1 \bigl[\bigl\{ m(\boldsymbol u,t) - \bar m(
\boldsymbol u)\bigr\}^2 + \sigma(\boldsymbol u,t)^2\bigr]
f(\boldsymbol u,t) \,dt.
\]
\end{theorem}

For any point $\boldsymbol u \in\mathbb{R}^d$ with $\bar
f(\boldsymbol
u) > 0$, we have
\[
E\bigl\{\tilde f_n(\boldsymbol u)\bigr\} = \bar f(\boldsymbol u) +
\frac{h_n^2 }{
2} \mathrm{tr} \biggl\{\int_0^1
H_f(\boldsymbol u,t) \bolds\kappa_S \,dt \biggr\} +
O\bigl(h_n^3\bigr)
\]
and
\begin{eqnarray*}
\frac{E\{\tilde T_n(\boldsymbol u)\} }{ E\{\tilde f_n(\boldsymbol u)\}
} & = & \bar m(\boldsymbol u) + \frac{h_n^2 }{2\bar f(\boldsymbol u)}
\mathrm{tr}
\biggl[\int_0^1 \bigl\{H_{mf}(
\boldsymbol u,t) - m(\boldsymbol u,t)H_f(\boldsymbol u,t)\bigr\}
\bolds\kappa_S \,dt \biggr]
\\
& &{} + O\bigl(h_n^3\bigr).
\end{eqnarray*}
Therefore, (\ref{eqnftilde}) and (\ref{eqnmuhat}) provide consistent
estimators of $\bar f$ and $\bar m$,\break (weighted) temporal averages of
the local density function $f$ and the regression function $m$,
respectively. For stationary processes, Theorem~\ref
{thmCLTftildemuhat} relates to traditional results on nonparametric
kernel estimators; see, for example, \citet{Robinson1983}, \citet
{Bosq1996} and \citet{Wu2005}. The AMSE optimal bandwidth for the
time-constant kernel estimators (\ref{eqnftilde}) and (\ref
{eqnmuhat}) satisfies $h_n \asymp n^{-1/(d+4)}$.

\subsection{Uniform bounds}\label{subsecunifbnds}
For stationary or independent observations, uniform bounds for kernel
estimators have been obtained by \citet{Peligrad1991}, \citet
{Andrews1995}, \citet{Bosq1996}, \citet{Masry1996}, \citet
{FanYao2003} and \citet{Hansen2008}, among others. \citet
{HallMullerWu2006} obtained a uniform bound for time-varying kernel
density estimators for independent observations, while \citet
{Vogt2012} considered kernel regression estimators under strong mixing
conditions. We shall here establish uniform bounds for the time-varying
kernel estimators (\ref{eqnfhat}) and (\ref{eqnmhat}) under the
locally strict stationarity framework (\ref{eqnframeworkx}). We need
the following assumptions:
\begin{longlist}[(A4)]
\item[(A4)] there exists a $q > 2$ such that $\Psi_{0,q} < \infty$ and
$\Psi_{m,q} = O(m^{-\alpha})$ for some $\alpha> 1/2 - 1/q$;
\item[(A5)] let $\mathscr X \subseteq\mathbb R^d$ be a compact set,
and assume $\inf_{t \in[0,1]} \inf_{\boldsymbol u \in\mathscr X}
f(\boldsymbol u,t) > 0$.
\end{longlist}

%
\begin{theorem}\label{thmunifbnds}
Assume \textup{(A1)}, \textup{(A3)--(A5)}, $b_n \to0$, $h_n \to0$ and $nb_nh_n^d \to
\infty$. \textup{(i)} If there exists $r > r' > 0$ such that $\sup_{t \in[0,1]}
\|\boldsymbol G(t;\bolds{\mathcal{ H}}_0) \|_r < \infty$ and
$n^{2/r'+2+d-q} b_n^{d-q} h_n^{d(d+q)} \to0$, then
\[
\sup_{t \in[0,1]} \sup_{\boldsymbol u \in\mathbb{R}^d} \bigl|\hat f_n(
\boldsymbol u,t) - E\bigl\{\hat f_n(\boldsymbol u,t)\bigr\}\bigr| =
O_p \biggl\{ \frac{(\log n)^{1/2} }{(nb_nh_n^d)^{1/2}} \biggr\}.
\]
\textup{(ii)} If $\eta_i \in\mathcal L^p$ for some $p > 2$, and $n^{2+d-q}
b_n^{d-q} h_n^{d(d+q)} \to0$, then
\[
\sup_{t \in[0,1]} \sup_{\boldsymbol u \in\mathscr X} \biggl\vert
\hat
m_n(\boldsymbol u,t) - \frac{E\{\hat T_n(\boldsymbol u,t)\} }{ E\{\hat
f_n(\boldsymbol u,t)\}} \biggr\vert=
O_p \biggl\{\frac{(\log n)^{1/2} }{
(nb_nh_n^d)^{1/2}} + \frac{n^{1/p}\log n }{ nb_nh_n^d} \biggr\}.
\]
\end{theorem}

If the bandwidths $b_n \asymp n^{-1/(d+5)}$ and $h_n \asymp
n^{-1/(d+5)}$ have the optimal AMSE rate, and $\eta_i \in\mathcal L^p$
for some $p > (d+5)/2$, then the bound in Theorem~\ref{thmunifbnds}(ii) can be simplified to $O_p\{(nb_nh_n^d)^{-1/2} (\log n)^{1/2}\}$.
Theorem~\ref{thmunifbndfmu} provides a uniform bound for (\ref
{eqnftilde}) and (\ref{eqnmuhat}).

%
\begin{theorem}\label{thmunifbndfmu}
Assume \textup{(A1)}, \textup{(A3)--(A5)}, $h_n \to0$ and $nh_n^d \to\infty$. \textup{(i)} If
there exists $r > r' > 0$ such that $\sup_{t \in[0,1]}\|\boldsymbol
G(t;\bolds{\mathcal{ H}}_0) \|_r < \infty$ and $n^{2/r'+2+d-q}\times\break
h_n^{d(d+q)} \to0$, then
\[
\sup_{\boldsymbol u \in\mathbb{R}^d} \bigl|\tilde f_n(\boldsymbol u) -
E\bigl\{
\tilde f_n(\boldsymbol u)\bigr\}\bigr| = O_p \biggl\{
\frac{(\log n)^{1/2} }{
(nh_n^d)^{1/2}} \biggr\}.
\]
\textup{(ii)} If $\eta_i \in\mathcal L^p$ for some $p > 2$, and $n^{2+d-q}
h_n^{d(d+q)} \to0$, then
\[
\sup_{\boldsymbol u \in\mathscr X} \biggl\vert\hat\mu
_n(\boldsymbol u)
- \frac{E\{\tilde T_n(\boldsymbol u)\} }{ E\{\tilde f_n(\boldsymbol
u)\}} \biggr\vert= O_p \biggl\{\frac{(\log n)^{1/2}
}{(nh_n^d)^{1/2}} +
\frac{n^{1/p}\log n }{ nh_n^d} \biggr\}.
\]
\end{theorem}

If the bandwidth $h_n \asymp n^{-1/(d+4)}$ is AMSE-optimal, and $\eta_i
\in\mathcal L^p$ for some $p > (d+4)/2$, then the bound in
Theorem~\ref
{thmunifbndfmu}(ii) can be simplified to\break  $O_p\{(nh_n^d)^{-1/2} (\log
n)^{1/2}\}$.

\subsection{Model selection}\label{subsecModelSelection}
Model $\mathrm I$ is quite general in the sense that it does not impose
any specific parametric form on the regression function and allows it
to change over time. However, in practice it is useful to check whether
model $\mathrm I$ can be reduced to its simpler special cases, namely
models $\mathrm{\II}$--$\mathrm{\IV}$. Model selection between models
$\mathrm{\II}$ and $\mathrm{\IV}$, or between models~$\mathrm{\III}$
and $\mathrm{\IV}$, has been studied in the literature mainly by means
of hypothesis testing; see references in Section~\ref{secintro}.
Nevertheless, less attention has been paid to distinguishing between
models $\mathrm{\II}$ and~$\mathrm{\III}$. We shall here propose an
information criterion that can consistently select the underlying true
model among candidate models $\mathrm I$--$\mathrm{\IV}$. Let
$\mathscr
T \subset(0, 1)$ be a compact set and $\mathscr I_n = \{i = 1, \ldots,
n \mid i/n \in\mathscr T\}$. We consider the restricted residual sum
of squares for model $\mathrm I$, which takes the form
\[
\textsc{rss}_n(\mathscr X,\mathscr T,\mathrm I) = \sum
_{i \in\mathscr
I_n} \bigl\{y_i - \hat m_n(
\boldsymbol x_i, i/n)\bigr\}^2 \indicator{\boldsymbol
x_i \in\mathscr X},
\]
where $\indicator{\cdot}$ is the indicator function. Similarly, we can
define $\textsc{rss}_n(\mathscr X,\mathscr T,\mathrm{\II})$,
$\textsc
{rss}_n(\mathscr X,\mathscr T,\mathrm{\III})$ and $\textsc
{rss}_n(\mathscr X,\mathscr T,\mathrm{\IV})$ for models $\mathrm{\II
}$--$\mathrm{\IV}$, respectively. For the simple linear regression model
$\mathrm{\IV}$, the parameter $\bolds\theta$ can be estimated by
the least squares estimate
%
%
\begin{equation}
\label{eqnthetahat} \hat{\bolds\theta}_n = \Biggl(\frac{1 }{ n}
\sum_{i=1}^n \boldsymbol x_i
\boldsymbol x_i^\top\Biggr)^{-1} \Biggl(
\frac{1 }{ n} \sum_{i=1}^n \boldsymbol
x_i y_i \Biggr).
\end{equation}
For the time-varying coefficient model $\mathrm{\III}$, let
$K_{T,b_n}(\cdot) = b_n^{-1} K_T(\cdot/b_n)$, and we can use the kernel
estimator of \citet{PriestleyChao1972},
%
%
\begin{equation}
\label{eqnbetahat} \hat{\bolds\beta}_n(t) = \Biggl\{
\frac{1 }{ n} \sum_{i=1}^n \boldsymbol
x_i \boldsymbol x_i^\top K_{T,b_n}(i/n-t)
\Biggr\}^{-1} \Biggl\{\frac{1 }{ n} \sum
_{i=1}^n \boldsymbol x_i
y_i K_{T,b_n}(i/n-t) \Biggr\}.
\end{equation}
For a candidate model $\varrho\in\{\mathrm I,\mathrm{\II},\mathrm
{\III
},\mathrm{\IV}\}$, we define the generalized information criterion
%
%
\begin{equation}
\label{eqnGIC} \textsc{gic}_{\mathscr X, \mathscr T}(\varrho) = \log
\bigl\{\textsc
{rss}_n(\mathscr X,\mathscr T,\varrho)/n\bigr\} + \tau_n
\textsc{df}(\varrho),
\end{equation}
where $\tau_n$ is a tuning parameter indicating the amount of
penalization and $\textsc{df}(\varrho)$ represents the model complexity
for model $\varrho\in\{\mathrm I,\mathrm{\II},\mathrm{\III
},\mathrm
{\IV}\}$ determined as follows. For the simple linear regression model
$\mathrm{\IV}$, following the convention we set the model complexity or
degree of freedom to be the number of potential predictors, namely
$\textsc{df} (\mathrm{\IV}) = d$. For the time-varying coefficient
model $\mathrm{\III}$, the effective number of parameters used in
kernel smoothing is $b_n^{-1}$ for each one of the $d$ predictors [see,
e.g., \citet{HurvichSimonoffTsai1998}], and thus we set $\textsc
{df}(\mathrm{\III}) = b_n^{-1} \,d$. Let $\textsc{iqr}_k$, $k =
1,\ldots
,d$, be the componentwise interquartile ranges of $(\boldsymbol x_i)$,
and motivated by the same spirit as in \citet
{HurvichSimonoffTsai1998}, we set $\textsc{df}(\mathrm{\II}) =
(h_n^{d})^{-1} \prod_{k=1}^d (2\textsc{iqr}_k)$ and $\textsc
{df}(\mathrm I) = (b_n h_n^d)^{-1} \prod_{k=1}^d (2\textsc{iqr}_k)$,
where $2\textsc{iqr} = 1$ for random variables having a uniform
distribution on $[0,1]$. The final model is selected by minimizing the
information criterion (\ref{eqnGIC}). We shall make the following assumption:
\begin{longlist}[(A6)]
\item[(A6)] eigenvalues of $\boldsymbol M(\boldsymbol G,t) = E \{
\boldsymbol G(t;\bolds{\mathcal{ H}}_0) \boldsymbol G(t;
\bolds{\mathcal{ H}}_0)^\top\}$ are bounded away from zero and infinity on $[0,1]$.
\end{longlist}

In order to establish the selection consistency of (\ref{eqnGIC}), in
addition to the results developed in Sections~\ref{subsecCLTs} and
\ref
{subsecunifbnds} regarding models $\mathrm I$ and $\mathrm{\II}$, we
need the following conditions on estimators (\ref{eqnthetahat}) and
(\ref{eqnbetahat}) for models $\mathrm{\IV}$ and $\mathrm{\III}$,
respectively:
\begin{longlist}[(P1)]
\item[(P1)] There exists a nonrandom sequence $\bolds\theta_n$
such that $\hat{\bolds\theta}_n - \bolds\theta_n =
O_p(n^{-1/2})$. If model $\mathrm{\IV}$ is correctly specified, then
$\bolds\theta_n$ can be replaced by the true value $\bolds\theta_0$.
\item[(P2)] There exists a sequence of nonrandom functions
$\bolds\beta_n\dvtx[0,1] \to\mathbb R^d$ such that
\[
\sup_{t \in\mathscr T} \bigl|\hat{\bolds\beta}_n(t) -
\bolds\beta_n(t)\bigr| = O_p(\phi_n),
\]
where $\phi_n = (nb_n)^{-1/2}(\log n)^{1/2} + b_n^2$. If model
$\mathrm
{\III}$ is correctly specified, then $\bolds\beta_n(\cdot)$ can
be replaced by the true coefficient function $\bolds\beta
_0(\cdot
)$ and
\[
\sup_{t \in\mathscr T} \biggl|\boldsymbol M(\boldsymbol G,t) \biggl\{\hat{
\bolds\beta}_n(t) - \bolds\beta_0(t) -
\frac{\kappa_T b_n^2
\bolds\beta''_0(t) }{2} \biggr\}
- \frac{1 }{ n} \sum_{i=1}^n
\boldsymbol x_i e_i K_{T,b_n} (i/n-t )\biggr| =
O_p\bigl(\phi_n^2\bigr),
\]
where $\boldsymbol x_i e_i \in\mathcal L^2$, $i = 1,\ldots,n$.
\end{longlist}

%
\begin{remark}Conditions (P1) and (P2) can be verified for locally stationary
processes with short-range dependence. For example, for the linear
regression model $\mathrm{\IV}$, by Lemma~5.1 of \citet{ZhangWu2012},
we have $\sum_{i=1}^n \{\boldsymbol x_i \boldsymbol x_i^\top-
E(\boldsymbol x_i \boldsymbol x_i^\top)\} = O_p(n^{1/2})$ and $\sum
_{i=1}^n \{\boldsymbol x_i y_i - E(\boldsymbol x_i y_i)\} =
O_p(n^{1/2})$. Hence we can use
\[
\bolds\theta_n = \Biggl\{\frac{1 }{ n} \sum
_{i=1}^n E\bigl(\boldsymbol x_i
\boldsymbol x_i^\top\bigr) \Biggr\}^{-1} \Biggl\{
\frac{1 }{ n} \sum_{i=1}^n E(
\boldsymbol x_i y_i) \Biggr\},
\]
which equals to $\bolds\theta_0$ if $y_i = \boldsymbol x_i^\top
\bolds\theta_0 + e_i$, $i = 1,\ldots,n$. This verifies condition
(P1). For the time-varying coefficient model $\mathrm{\III}$, by
Lemma~5.3 of \citet{ZhangWu2012}, we have $\sup_{t \in\mathscr T} |n^{-1}
\sum_{i=1}^n \{\boldsymbol x_i \boldsymbol x_i^\top- E(\boldsymbol x_i
\boldsymbol x_i^\top)\} K_{T,b_n}(i/n-t)| = O_p(\phi_n)$ and $\sup_{t
\in\mathscr T} |n^{-1} \sum_{i=1}^n \{\boldsymbol x_i y_i -
E(\boldsymbol x_i y_i)\} K_{T,b_n}(i/n-t)| = O_p(\phi_n)$. Hence we
can use
\[
\bolds\beta_n(t) = \Biggl\{\frac{1 }{ n} \sum
_{i=1}^n E\bigl(\boldsymbol x_i
\boldsymbol x_i^\top\bigr) K_{T,b_n}(i/n-t) \Biggr
\}^{-1} \Biggl\{\frac{1 }{
n} \sum_{i=1}^n
E(\boldsymbol x_i y_i) K_{T,b_n}(i/n-t) \Biggr
\},
\]
and condition (P2) follows by the proof of Theorem~3 in \citet{ZhouWu2010}.
\end{remark}

Recall that the AMSE optimal bandwidths satisfy $b_n(\mathrm I) \asymp
n^{-1/(d+5)}$ and $h_n(\mathrm I) \asymp n^{-1/(d+5)}$ for model
$\mathrm I$, $h_n(\mathrm{\II}) \asymp n^{-1/(d+4)}$ for model
$\mathrm
{\II}$ and $b_n(\mathrm{\III}) \asymp n^{-1/5}$ for model $\mathrm
{\III
}$. Theorem~\ref{thmGIC} provides the selection consistency of the
information criterion (\ref{eqnGIC}), where the true model is denoted
by $\varrho_0$.

%
\begin{theorem}\label{thmGIC}
Assume \textup{(A1)}, \textup{(A3)--(A6)} with $q > (3d+5)/(d+2)$, \textup{(P1)},
\textup{(P2)}, $\eta_i \in
\mathcal L^p$ for some $p > (d+5)/2$, $i = 1,\ldots,n$, and bandwidths
with optimal AMSE rates are used for models $\mathrm I$--$\mathrm{\III
}$. If
\[
\tau_n n^{(d+1)/(d+5)} \to0, \qquad\tau_n n^{(d+3)/(d+4)}
\to\infty,
\]
then for any $\varrho_1 \in\{\mathrm I,\mathrm{\II},\mathrm{\III
},\mathrm{\IV}\}$ and $\varrho_1 \neq\varrho_0$, we have
\[
\mathrm{pr}\bigl\{\textsc{gic}_{\mathscr X, \mathscr T}(\varrho_0)
< \textsc
{gic}_{\mathscr X, \mathscr T}(\varrho_1)\bigr\} \to1.
\]
\end{theorem}

\subsection{Extensions}\label{subsecextensions}
Recall that in Theorems \ref{thmCLTfhatmhat}--\ref{thmGIC} error
process (\ref{eqndefeiiideta}) has i.i.d. $\eta_i$, which are also
independent of $(\boldsymbol x_j)$. In Section~\ref
{subsubsecextensioncorrelation} we allow serially correlated~$\eta_i$.
Section~\ref{subsubsecextensionautoregression} concerns time-varying
autoregressive processes in which $(\eta_i)$ and $(\boldsymbol x_j)$
are naturally dependent.

\subsubsection{Models with serially correlated errors}\label
{subsubsecextensioncorrelation}
To allow errors with serial correlation, similarly to (\ref
{eqnframeworkx}) we assume that
%
%
\begin{equation}
\label{eqnetai} \eta_i = L(i/n; \mathcal J_i),
\end{equation}
where $\mathcal J_i = (\ldots, \zeta_{i-1}, \zeta_i)$ with $\zeta_k$,
$k \in\mathbb Z$, being i.i.d. random variables and independent of
$\bolds\xi_j$, $j \in\mathbb Z$. Therefore, $(\eta_i)$ is a
dependent nonstationary process that is independent of $(\boldsymbol
x_j)$, and the error process $e_i = \sigma(\boldsymbol x_i,i/n)\eta_i$
can exhibit both serial correlation and heteroscedasticity; see \citet
{Robinson1983}, \citeauthor{OrbeFerreiraRodriguezPoo2005}
(\citeyear{OrbeFerreiraRodriguezPoo2005,OrbeFerreiraRodriguezPoo2006}) and
references therein for similar error structures. Let $\zeta_i', \zeta
_j, i,j \in\mathbb Z$, be i.i.d. and $\mathcal J_k' = (\ldots,\zeta
_{-1},\zeta_0',\zeta_1,\break \ldots,\zeta_k)$. Assume $c_{L,q} = \sup_{t
\in
[0,1]} \|L(t;\mathcal J_0)\|_q < \infty$, and define the functional
dependence measure
\[
\nu_{k,q} = \sup_{t \in[0,1]} \bigl\|L(t; \mathcal
J_k) - L\bigl(t; \mathcal J_k'\bigr)
\bigr\|_q.
\]
The following theorem states that the results presented in
Sections~\ref{subsecCLTs}--\ref{subsecModelSelection} will continue
to hold
(except for a difference of $\log n$ on the uniform bounds) if the
process $(\eta_i)$ in (\ref{eqnetai}) satisfies the geometric moment
contraction (GMC) condition [\citet{ShaoWu2007}]. The proof is
available in the supplementary material [\citet{ZhangWu2014}].

%
\begin{theorem}\label{thmextensioncorrelation}
Assume that the process $(\eta_i)$ in (\ref{eqnetai}) satisfies $\nu
_{k,4} = O(\rho^k)$ for some $0 < \rho< 1$. Then the results of
Theorems \ref{thmCLTfhatmhat}--\ref{thmGIC} will continue to hold
except that the uniform bounds in Theorems~\ref{thmunifbnds}\textup{(ii)}
and~\ref{thmunifbndfmu}\textup{(ii)} will be multiplied by a factor of
$\log n$.
\end{theorem}

\subsubsection{Time-varying nonlinear autoregressive models}
\label{subsubsecextensionautoregression}
In this section we shall consider the autoregressive version of (\ref
{eqnformmxt}),
%
%
\begin{eqnarray}
\label{eqnS10754a} y_i = m(\boldsymbol x_i,i/n) + \sigma(
\boldsymbol x_i,i/n)\eta_i,
\nonumber
\\[-8pt]
\\[-8pt]
\eqntext{ \boldsymbol x_i
= (y_{i-1},\ldots,y_{i-d})^\top, i = 1,\ldots,n,}
\end{eqnarray}
where $\eta_i$ are i.i.d. random variables with $E(\eta_i) = 0$ and
$E(\eta_i^2) = 1$. We can view~(\ref{eqnS10754a}) as a time-varying or
locally stationary autoregressive process, and the corresponding shift
processes $\mathcal F_k = (\ldots,\eta_{k-1},\eta_k)$ and $\mathcal H_k
= \mathcal F_{k-1}$. We shall here present analogous versions of
Theorems \ref{thmCLTfhatmhat}--\ref{thmGIC}. Note that in this case
${\boldsymbol x_i}$ cannot be written in the form of (\ref
{eqnframeworkx}). However, Proposition~\ref{propS10828a} implies that
it can be well approximated by a process in the form of (\ref
{eqnframeworkx}). For each $t \in[0, 1]$, we define the process $\{
y_i(t)\}_{i \in\mathbb Z}$ by
%
%
\begin{eqnarray}
\label{eqnS10755a} y_i(t) &=& m\bigl\{\boldsymbol x_i(t),t
\bigr\} + \sigma\bigl\{\boldsymbol x_i(t),t\bigr\}\eta
_i,
\nonumber
\\[-8pt]
\\[-8pt]
\nonumber
\boldsymbol x_i(t)& =& \bigl\{y_{i-1}(t),
\ldots,y_{i-d}(t)\bigr\}^\top.
\end{eqnarray}

%
\begin{lemma}\label{lemS10343}
Assume that there exist constants $a_1, \ldots, a_d \ge0$ with\break $\sum
_{j=1}^d a_j < 1$, such that, for all $\boldsymbol x = (x_1,\ldots
,x_d)^\top$ and $\boldsymbol x' = (x'_1,\ldots,x'_d)^\top$,
%
%
\begin{eqnarray}
\label{eqS11719p} &&\sup_{0 \le t \le1} \bigl\| \bigl[m({\boldsymbol x},t) +
\sigma({\boldsymbol x},t) \eta_i\bigr] - \bigl[m\bigl({\boldsymbol
x}',t\bigr) + \sigma\bigl({\boldsymbol x}',t\bigr)
\eta_i\bigr]\bigr \|_p
\nonumber
\\[-8pt]
\\[-8pt]
\nonumber
&&\qquad\le\sum_{j=1}^d
a_j \bigl|x_j - x_j'\bigr|.
\end{eqnarray}
Then \textup{(i)} the recursion (\ref{eqnS10755a}) has a stationary solution of
the form $y_i(t) = g(t;\mathcal F_i)$ which satisfies the geometric
moment contraction (GMC) property: for some $\rho\in(0, 1)$,
\[
\sup_{0 \le t \le1} \delta_i(t) = O\bigl(
\rho^i\bigr), \qquad\delta_i(t) = \bigl\| g(t;\mathcal
F_i) - g\bigl(t;\mathcal F'_i\bigr)
\bigr\|_p.
\]
\textup{(ii)} If in (\ref{eqnS10754a}) the initial values $(y_0, y_{-1}, \ldots
, y_{1-d}) = {\boldsymbol x_1}(0)$, then $y_i$ can be written in the
form $g_i (\mathcal F_i)$, where $g_i(\cdot)$ is a measurable function,
and it also satisfies the GMC property
%
%
\begin{equation}\qquad
\label{eqS11707p} \sup_{i \le n} \bigl\| y_i-
g_i\bigl(\ldots, \eta_{i-k-2}, \eta_{i-k-1}, \eta
_{i-k}', \eta_{i-k+1}, \ldots,
\eta_i \bigr) \bigr\|_p = O\bigl(\rho^k\bigr).
\end{equation}
\end{lemma}

Lemma~\ref{lemS10343}(i) concerns the stationarity of the process $\{
y_i(t)\}_{i \in\mathbb Z}$, which follows from Theorem~5.1 of \citet
{ShaoWu2007}. For (ii), denote by $\theta_k^\dag$ the left-hand side
of (\ref{eqS11707p}). Then by (\ref{eqS11719p}), $\theta_k^\dag$
satisfies $\theta_k^\dag\leq\sum_{j=1}^d a_j \theta_{k-j}^\dag$,
implying (\ref{eqS11707p}) via recursion.

For presentational simplicity suppose we observe $y_{1-d}, y_{2-d},
\ldots, y_n$ from model (\ref{eqnS10754a}) with the initial values
$(y_0, y_{-1}, \ldots, y_{1-d}) = {\boldsymbol x_1}(0)$. Estimates
(\ref
{eqnfhat}) and (\ref{eqnmhat}) can be computed in the same way.
Proposition~\ref{propS10828a} implies that, for $i$ such that $i/n
\approx u$, the process $({\boldsymbol x}_i)_i$ can be approximated by
the stationary process $\{\boldsymbol x_i(u)\}_i$, thus suggesting
local strictly stationarity. The proof is available in the supplementary material
[\citet{ZhangWu2014}].

%
\begin{proposition}
\label{propS10828a}
Let $G_\eta({\boldsymbol x}, t) = m({\boldsymbol x}, t) + \sigma
({\boldsymbol x}, t) \eta$ and $\dot G_\eta({\boldsymbol x}, t) =\break
\partial G_\eta({\boldsymbol x}, t) / \partial t$. Assume (\ref
{eqS11719p}) and
\[
\sup_{0 \le t \le1} \sup_{0 \le u \le1} \bigl\| \dot
G_{\eta_i}\bigl\{ {\boldsymbol x}_i(u), t\bigr\}
\bigr\|_p < \infty.
\]
Then $\| {\boldsymbol x}_i - {\boldsymbol x}_i(u) \|_p = O(n^{-1} + |u-i/n|)$.
\end{proposition}

Let $f(\boldsymbol u, t)$ be the density of ${\boldsymbol x_i}(t) = \{
y_{i-1}(t), \ldots, y_{i-d}(t)\}$ and $f_\eta$ be the density of
$\eta
_i$. Theorem~\ref{thmextensionautoregression} serves as an analogous
version of Theorems \ref{thmCLTfhatmhat}--\ref{thmunifbndfmu}, and
the proof is available in the supplementary material [\citet{ZhangWu2014}].

%
\begin{theorem}\label{thmextensionautoregression}
Assume \textup{(A1)}, \textup{(A5)} and $\sup_w \{f_\eta(w) + |f_\eta'(w)|\} < \infty$.
Let the conditions in Lemma~\ref{lemS10343} and Proposition~\ref
{propS10828a} be satisfied. Then under respective conditions in
Theorems \ref{thmCLTfhatmhat}--\ref{thmGIC}, the corresponding
conclusions also hold, respectively.
\end{theorem}

\section{Numerical implementation}\label{secnumerexp}
\subsection{Bandwidth and tuning parameter selection}\label{subsecbandwidth}
Selecting bandwidths that optimize the performance of (\ref{eqnGIC})
can be quite nontrivial, and in our case, it is further complicated by
the presence of dependence and nonstationarity. Assuming independence,
the problem of bandwidth selection has been considered for model~$\mathrm{\II}$ by \citet{HardleMarron1985}, \citet
{HardleHallMarron1988}, \citet{ParkMarron1990}, \citet
{RuppertSheatherWand1995}, \citet{WandJones1995}, \citet{Xia1998}
and \citet{GaoGijbels2008}, among others. \citet
{HooverRiceWuYang1998}, \citet{FanZhangJT2000} and \citet
{RamsaySilverman2005} considered the problem for model $\mathrm{\III
}$ for longitudinal data, where multiple independent realizations are
available. For the time-varying kernel density estimator (\ref
{eqnfhat}) with independent observations, \citet{HallMullerWu2006}
coupled the selection of spatial and temporal bandwidths and adopted
the least squares cross validation [\citet{Silverman1986}].
Nevertheless, bandwidths selectors derived under independence can break
down for dependent data [\citet{Wang1998} and \citet
{OpsomerWangYang2001}]. We propose using the AMSE optimal bandwidths
$b_n(\mathrm I) = c_b(\mathrm I) n^{-1/(d+5)}$ and $h_n(\mathrm I) =
c_h(\mathrm I) n^{-1/(d+5)}$ for model $\mathrm I$, $h_n(\mathrm{\II})
= c_h(\mathrm{\II}) n^{-1/(d+4)}$ for model $\mathrm{\II}$ and
$b_n(\mathrm{\III}) = c_b(\mathrm{\III}) n^{-1/5}$ for model
$\mathrm
{\III}$, where $0 < c_b(\mathrm I), c_h(\mathrm I), c_h(\mathrm{\II}),
c_b(\mathrm{\III}) < \infty$ are constants. Due to the presence of both
dependence and nonstationarity, estimation of these constants is
difficult. Throughout this section, as a rule of thumb, we use
$c_b(\mathrm I) = c_b(\mathrm{\III}) = 1/2$ and $c_h(\mathrm I) =
c_h(\mathrm{\II}) = \prod_{k=1}^d \textsc{iqr}_k$. Our numerical
examples suggest that these simple choices have a reasonably good performance.

We shall here discuss the choice of the tuning parameter $\tau_n$ that
controls the amount of penalization on models complexities. The problem
has been extensively studied for the linear model $\mathrm{\IV}$ by
\citet{Akaike1973}, \citet{Mallows1973}, \citet{Schwarz1978}, \citet
{Shao1997} and \citet{Yang2005} among others. For the generalized
information criterion (\ref{eqnGIC}), given conditions in Theorem~\ref
{thmGIC}, one can choose $\tau_n = c n^{-(d+3)/(d+4)} \log n$, where
$c > 0$ is a constant, which satisfies all the required conditions and
thus guarantees the selection consistency. Note that the choice of $c$
does not affect the asymptotic result, namely the proposed method will
select the true model for any given $c > 0$ as long as the sample size
is large enough; see Theorem~\ref{thmGIC}. Therefore, one can simply
use $c = 1$ to devise a consistent model selection procedure. As an
alternative, following \citet{FanLi2001} and \citet
{TibshiraniTibshirani2009}, we shall here consider a data-driven
selector based on the $K$-fold cross-validation (CV). In particular, we
first split the data into $K$ parts, denoted by $\mathcal D_1,\ldots
,\mathcal D_K$, then for each $k = 1,\ldots,K$, we remove the $k$th
part from the data and use the information criterion (\ref{eqnGIC}) to
select the model, based on which predictions can be made for the
removed part and are denoted by $\hat y_i^{-k}(c)$, $i \in\mathcal
D_k$. The selected value $\hat c$ is obtained by minimizing the
cross-validation criterion
\[
\textsc{cv}(c) = \sum_{k=1}^K \sum
_{i \in\mathcal D_k} \bigl\{y_i - \hat
y_i^{-k}(c)\bigr\}^2.
\]
It can be seen from the simulation results in Section~\ref
{subsecsimulation} that this CV-based tuning parameter selector
performs reasonably well.

\subsection{Simulation results}\label{subsecsimulation}
We shall in this section carry out a simulation study to examine the
finite-sample performance of the generalized information criterion~(\ref
{eqnGIC}). Let $d = 1$ and $\xi_i$, $i \in\mathbb Z$ and $\eta_j$, $j
\in\mathbb Z$ be i.i.d. standard normal, $a(t) = (t-1/2)^2$, $t \in
[0,1]$ and $G(t; \mathcal H_k) = \xi_k + \sum_{l=1}^\infty a(t)^l \xi
_{k-l}$, $k \in\mathbb Z$, $t \in[0,1]$. For the regressor and error
processes with $x_i = G(i/n; \mathcal H_i)$ and $e_i = \sigma
(x_i,i/n)\eta_i$, $i = 1,\ldots,n$, we consider model (\ref
{eqnformmxt}) with the following four specifications:
\begin{longlist}[(a)]
\item[(a)] $m(x,t) = 2.5\sin(2\pi t)\cos(\pi x)$ and $\sigma(x,t) =
\varphi|tx|/2$;
\item[(b)] $m(x,t) = \exp(x)$ and $\sigma(x,t) = \varphi t \exp(x/3)$;
\item[(c)] $m(x,t) = 5t + 4\cos(2 \pi t)x$ and $\sigma(x,t) =
\varphi
\exp(tx/2)$;
\item[(d)] $m(x,t) = 2 + 3x$ and $\sigma(x,t) = \varphi|x/3+t|$,
\end{longlist}
where $\varphi> 0$ is a constant indicating the noise level. Cases
(a)--(d) correspond to models $\mathrm I$--$\mathrm{\IV}$, respectively,
and their signal-to-noise ratios (SNRs) are roughly of the same order
given the same $\varphi$. The Epanechnikov kernel $K(v) = 3(1-v^2)/4$,
$v \in[-1,1]$, is used hereafter for both the spatial and temporal
kernel functions. Let $\mathscr X = [-2,2]$ and $\mathscr T =
[0.2,0.8]$. The tuning parameter is selected by using the tenfold
CV-based method described in Section~\ref{subsecbandwidth}. The
results are summarized in Table~\ref{tabGICsimulation} for different
noise levels $\varphi\in\{1,2,3\}$ and sample sizes $n = 2^k \times
250$, $0 \leq k \leq3$. For each configuration, the results are based
on 1000 simulated realizations of models (a)--(d).

\begin{sidewaystable}
\tabcolsep=0pt
\tablewidth=\textwidth
\caption{Proportions of selecting models $\mathrm I$--$\mathrm{\IV}$
for different combinations of noise levels $\varphi$, sample sizes $n$
and model specifications \textup{(a)--(d)} with 1000 replications for each
configuration. Medians of the SNR are also reported, where for each
realization $y_i = m_i(x_i) + e_i$, $i = 1,\ldots,n$, the SNR is
defined as $\{\sum_{i=1}^n m_i(x_i)^2 / \sum_{i=1}^n e_i^2\}
^{1/2}$}\label{tabGICsimulation}
\begin{tabular*}{\textwidth}{@{\extracolsep{\fill}}lcccccccccccccccccc@{}}
\hline
& & \multicolumn{5}{c}{$\bolds{\varphi= 1}$} & & \multicolumn{5}{c}{$\bolds{\varphi = 2}$} & & \multicolumn{5}{c@{}}{$\bolds{\varphi= 3}$} \\[-6pt]
& & \multicolumn{5}{c}{\hrulefill} & & \multicolumn{5}{c}{\hrulefill} & &
\multicolumn{5}{c@{}}{\hrulefill} \\
& & & \multicolumn{4}{c}{\textbf{Selected model}} & & & \multicolumn
{4}{c}{\textbf{Selected model}} & & & \multicolumn{4}{c@{}}{\textbf{Selected model}} \\[-6pt]
& & & \multicolumn{4}{c}{\hrulefill} & & & \multicolumn
{4}{c}{\hrulefill} & & & \multicolumn{4}{c@{}}{\hrulefill} \\
$\bolds{n}$ & \textbf{Case} & \textbf{SNR} & $ \mathbf I$ & $\mathbf{II}$ & $\mathbf{III}$ &
$\mathbf{IV}$ & & \textbf{SNR} & $\mathbf I$& $\mathbf{II}$ & $\mathbf{III}$
& $\mathbf{IV}$ & & \textbf{SNR} & $\mathbf I$ & $\mathbf{II}$ & $\mathbf{III}$ & \multicolumn{1}{c@{}}{$\mathbf{IV}$} \\
\hline
\phantom{0}250 & (a) & 4.36 & 0.967 & 0.000 & 0.000 & 0.033 & & 2.16 & 0.920 &
0.000 & 0.000 & 0.080 & & 1.45 & 0.840 & 0.000 & 0.000 & 0.160 \\
& (b) & 4.09 & 0.116 & 0.882 & 0.000 & 0.002 & & 2.04 & 0.119 & 0.857
& 0.000 & 0.024 & & 1.36 & 0.132 & 0.784 & 0.002 & 0.082 \\
& (c) & 3.73 & 0.016 & 0.000 & 0.984 & 0.000 & & 1.86 & 0.032 & 0.000
& 0.968 & 0.000 & & 1.24 & 0.032 & 0.000 & 0.968 & 0.000 \\
& (d) & 5.44 & 0.017 & 0.043 & 0.005 & 0.935 & & 2.72 & 0.014 & 0.040
& 0.001 & 0.945 & & 1.82 & 0.024 & 0.040 & 0.003 & 0.933 \\[3pt]
\phantom{0}500 & (a) & 4.29 & 0.985 & 0.000 & 0.000 & 0.015 & & 2.15 & 0.945 &
0.000 & 0.000 & 0.055 & & 1.44 & 0.896 & 0.000 & 0.000 & 0.104 \\
& (b) & 4.17 & 0.044 & 0.949 & 0.000 & 0.008 & & 2.08 & 0.058 & 0.906
& 0.000 & 0.036 & & 1.40 & 0.037 & 0.926 & 0.000 & 0.037 \\
& (c) & 3.71 & 0.001 & 0.000 & 0.999 & 0.000 & & 1.86 & 0.008 & 0.000
& 0.992 & 0.000 & & 1.24 & 0.012 & 0.000 & 0.988 & 0.000 \\
& (d) & 5.42 & 0.007 & 0.037 & 0.000 & 0.956 & & 2.71 & 0.012 & 0.042
& 0.001 & 0.945 & & 1.81 & 0.005 & 0.026 & 0.006 & 0.963 \\[3pt]
1000 & (a) & 4.29 & 0.994 & 0.000 & 0.000 & 0.006 & & 2.15 & 0.970 &
0.000 & 0.000 & 0.030 & & 1.44 & 0.921 & 0.000 & 0.000 & 0.079 \\
& (b) & 4.17 & 0.004 & 0.992 & 0.000 & 0.004 & & 2.08 & 0.005 & 0.975
& 0.000 & 0.020 & & 1.40 & 0.015 & 0.957 & 0.000 & 0.028 \\
& (c) & 3.71 & 0.000 & 0.000 & 1.000 & 0.000 & & 1.86 & 0.001 & 0.000
& 0.999 & 0.000 & & 1.24 & 0.004 & 0.000 & 0.996 & 0.000 \\
& (d) & 5.42 & 0.001 & 0.028 & 0.002 & 0.969 & & 2.71 & 0.002 & 0.024
& 0.003 & 0.971 & & 1.81 & 0.001 & 0.025 & 0.002 & 0.972 \\[3pt]
2000 & (a) & 4.29 & 0.999 & 0.000 & 0.000 & 0.001 & & 2.15 & 0.979 &
0.000 & 0.000 & 0.021 & & 1.44 & 0.948 & 0.000 & 0.000 & 0.052 \\
& (b) & 4.17 & 0.000 & 0.997 & 0.000 & 0.003 & & 2.08 & 0.000 & 0.982
& 0.000 & 0.018 & & 1.40 & 0.000 & 0.965 & 0.000 & 0.035 \\
& (c) & 3.71 & 0.000 & 0.000 & 1.000 & 0.000 & & 1.86 & 0.000 & 0.000
& 1.000 & 0.000 & & 1.24 & 0.001 & 0.000 & 0.999 & 0.000 \\
& (d) & 5.42 & 0.000 & 0.014 & 0.001 & 0.985 & & 2.71 & 0.000 & 0.014
& 0.001 & 0.985 & & 1.81 & 0.000 & 0.014 & 0.000 & 0.986 \\
\hline
\end{tabular*}
\end{sidewaystable}

It can be seen from Table~\ref{tabGICsimulation} that the proposed
model selection procedure performs reasonably well as it has very high
empirical probabilities of identifying the true model, even when the
sample size is moderate to small. For example, if the sample size $n =
250$, which is usually considered to be small for conducting
time-varying nonparametric inference, and the data are generated by
model (a) with $\varphi= 1$, then 967 out of 1000 realizations are
correctly identified as the time-varying nonparametric regression model
$\mathrm I$, while 33 out of 1000 realizations are under-fitted as the
simple linear regression model $\mathrm{\IV}$. Hence, for each
combination of $n$ and $\varphi$, in the ideal case, we expect the
block to have unit diagonal components and zero off-diagonal
components. For each configuration, medians of the SNR are also
reported, where for each realization $y_i = m_i(x_i) + e_i$, $i =
1,\ldots,n$, the SNR is defined as $\{\sum_{i=1}^n m_i(x_i)^2 / \sum
_{i=1}^n e_i^2\}^{1/2}$. It can be seen that the proposed model
selection procedure with the CV-based tuning parameter selector has a
reasonably robust performance with respect to the noise level, and the
performance improves quickly if we increase the sample size. Note that
a sample size of 1000 is considered to be reasonable if one would like
to conduct time-varying nonparametric inference.

\subsection{Application on modeling interest rates}
\label{subsecinterestrates} 
\begin{figure}

\includegraphics{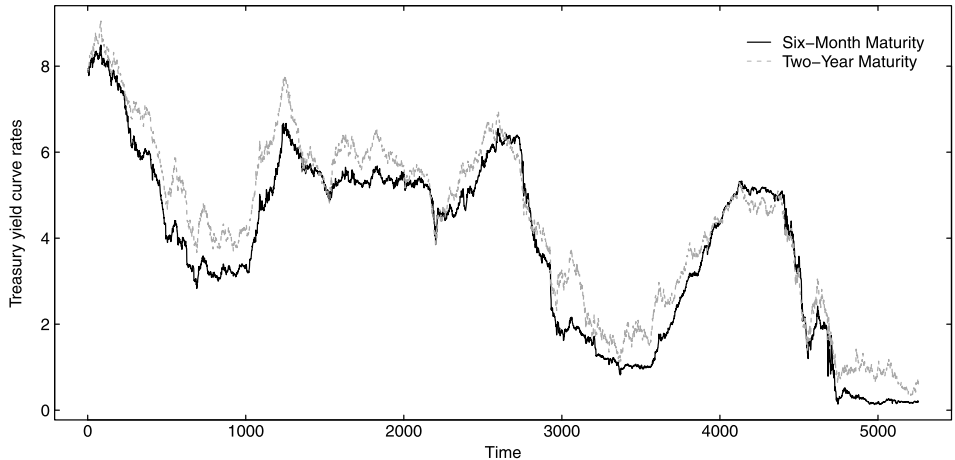}

\caption{Time series plots for the U.S. daily treasury yield curve
rates with six-month (solid black) and two-year (dashed grey)
maturities.}\label{figtreasuryyieldrates}
\end{figure}

Modeling interest rates is an
important problem in finance. In \citet{BlackScholes1973} and
\citet{Merton1974} interest rates were assumed to be constants. A
popular model is the time-homogeneous diffusion process
(\ref{eqnSDEx}) with linear drift function; see, for example,
\citet{Vasicek1977}, \citet{Courtadon1982},
\citet{CoxIngersollRoss1985} and
\citet{ChanKarolyiLongstaffSanders1992}. Its discretized
version is given by model $\mathrm{\IV}$.
\citet{AitSahalia1996}, \citet{Stanton1997} and
\citet{LiuWu2010} considered model (\ref{eqnSDEx}) with
nonlinear drift function, which relates to model $\mathrm{\II}$.
We consider the daily U.S. treasury yield curve rates with
six-month and two-year maturities during 01/02/1990--12/31/2010.
The data can be obtained from the U.S. Department of the Treasury
website at \url{http://www.treasury.gov/}. Both series contain $n
= 5256$ daily rates, and their time series plots are shown in
Figure~\ref{figtreasuryyieldrates}.

%

We shall here model the data by the time-varying diffusion process
(\ref
{eqnSDExt}), and apply the proposed model selection procedure to
determine the forms of the drift functions. Let $x_i = r_{t_i}$ be the
observation at day $i$. Since a year has 250 transaction days, $\Delta
= t_i - t_{i-1} = 1/250$.\vadjust{\goodbreak} Following \citet{LiuWu2010}, we consider
the following discretized version of (\ref{eqnSDExt}):
%
%
\begin{eqnarray}
\label{eqnSDExtdiscrete} y_i = r_{t_{i+1}} - r_{t_i} =
\mu(x_i,i/n)\Delta+ \sigma(x_i,i/n)
\Delta^{1/2} \eta_i
\nonumber
\\[-8pt]
\\[-8pt]
\eqntext{\mbox{where } \displaystyle\eta_i =
\frac{{\B_{t_{i+1}} - \B_{t_i}} }{\Delta^{1/2}}.}
\end{eqnarray}
Note that $\eta_i$ are i.i.d. $N\{0, 1\}$ random variables. We shall
here write\break $\mu(x_i, i/n) \Delta$ and $\sigma(x_i,i/n) \Delta
^{1/2}$ in
(\ref{eqnSDExtdiscrete}) as $m(x_i,i/n)$ and $\sigma(x_i,i/n)$ in the
sequel. Then specifications of \citet{Vasicek1977} and \citet
{LiuWu2010} become models $\mathrm{\IV}$ and $\mathrm{\II}$, respectively.

For the treasury yield curve rates with six-month maturity, let
$\mathscr T = [0.2,0.8]$, and $\mathscr X = [0.18,7.89]$ which includes
95.5\% of the daily rates $x_i$. The selected bandwidths and tuning
parameter are $b_n(\mathrm I) = 0.12$, $h_n(\mathrm I) = 0.82$,
$h_n(\mathrm{\II}) = 0.62$, $b_n(\mathrm{\III}) = 0.09$ and $\hat
\tau
_n = 0.00090$. The results are summarized in Table~\ref
{tabTreasuryYieldRates}. Hence, the time-varying coefficient model
$\mathrm{\III}$ is selected, and we conclude that the treasury yield
curve rates with six-month maturity should be modeled by (\ref
{eqnSDExt}) with $\mu(r_t,t) = \beta_0(t) + \beta_1(t) r_t$ for some
smoothly varying functions $\beta_0(\cdot)$ and $\beta_1(\cdot)$, which
serves as a time-varying version of \citet{ChanKarolyiLongstaffSanders1992}.

%
\begin{table}
\caption{Results of the model selection procedure based on the
generalized information criterion (\protect\ref{eqnGIC}) for treasury yield
rates with six-month and two-year maturity periods}\label{tabTreasuryYieldRates}
\begin{tabular*}{\textwidth}{@{\extracolsep{\fill}}ld{2.3}d{2.2}d{2.3}d{2.3}d{2.2}d{2.3}@{}}
\hline
& \multicolumn{3}{c}{\textbf{Six-month maturity}} & \multicolumn{3}{c}{\textbf{Two-year maturity}} \\[-6pt]
& \multicolumn{3}{c}{\hrulefill} & \multicolumn{3}{c@{}}{\hrulefill} \\
\textbf{Model} & \multicolumn{1}{c}{$\bolds{\log(\textsc{\textbf{rss}}/n)}$} & \multicolumn{1}{c}{$\textsc{\textbf{df}}$} & \multicolumn{1}{c}{$\textsc{\textbf{gic}}$} &
\multicolumn{1}{c}{$\bolds{\log(\textsc{\textbf{rss}}/n)}$} & \multicolumn{1}{c}{$\textsc{\textbf{df}}$} & \multicolumn{1}{c@{}}{$\textsc{\textbf{gic}}$} \\
\hline
$\mathrm I$ & -6.853 & 69.54 & -6.790 & -6.126 & 69.54 & -6.063 \\
$\mathrm{\II}$ & -6.824 & 11.10 & -6.814 & -6.114 & 11.10 & -6.104
\\
$\mathrm{\III}$ & -6.851 & 22.19 & -6.831 & -6.126 & 22.19 & -6.106
\\
$\mathrm{\IV}$ & -6.822 & 2.00 & -6.820 & -6.113 & 2.00 & -6.111 \\
\hline
\end{tabular*}
\end{table}

We then consider the treasury yield curve rates with two-year maturity.
Let $\mathscr T = [0.2,0.8]$ and $\mathscr X = [0.67,8.16]$ which
includes 95.1\% of the daily rates $x_i$. The selected bandwidths and
tuning parameter are $b_n(\mathrm I) = 0.12$, $h_n(\mathrm I) = 0.75$,
$h_n(\mathrm{\II}) = 0.56$, $b_n(\mathrm{\III}) = 0.09$ and $\hat
\tau
_n = 0.00090$. Based on Table~\ref{tabTreasuryYieldRates}, the linear
regression model $\mathrm{\IV}$ is selected. In comparison with the
results with six-month maturity, our analysis suggests that treasury
yield rates with longer maturity are more stable over time.

\section{Conclusion}\label{secconclusion}
The paper considers a time-varying nonparametric regression model,
namely model $\mathrm I$, which is able to capture time-varying and
nonlinear relationships between the response variable and the
explanatory variables. It includes the popular nonparametric regression
model $\mathrm{\II}$ and time-varying coefficient model~$\mathrm
{\III}$
as special cases, and all of them are generalizations of the simple
linear regression model $\mathrm{\IV}$. In comparison with existing
results, the current paper makes two major contributions. First, we
develop an asymptotic theory on nonparametric estimation of the
time-varying regression model (\ref{eqnformmxt}) under the new
framework of \citet{DraghicescuGuillasWu2009}. Compared with the
classical strong mixing conditions as used by \citet{Vogt2012}, the
current framework is convenient to work with and often leads to optimal
asymptotic results. In the proof, we use both the martingale
decomposition and the $m$-dependence approximation techniques to obtain
sharp results. Second, although the time-varying regression model
$\mathrm I$ is quite general by allowing a time-varying nonlinear
relationship between the response variable and the explanatory
variables, it can be useful in practice to check whether it can be
reduced to its simpler special cases, namely models $\mathrm{\II
}$--$\mathrm{\IV}$ which have been extensively used in the literature.
However, existing results on model selection usually focused on
distinguishing between models $\mathrm{\II}$ and $\mathrm{\IV}$ and
between models $\mathrm{\III}$ and $\mathrm{\IV}$, and much less
attention has been paid to distinguishing between models $\mathrm{\II}$
and $\mathrm{\III}$. Note that models $\mathrm{\II}$ and $\mathrm
{\III
}$ are both generalizations of the simple linear regression model
$\mathrm{\IV}$ but in completely different aspects, and therefore it is
desirable if we can have a statistically valid method to decide which
generalization (or the more general model $\mathrm I$) should be used
for a given data set. The current paper fills this gap by proposing an
information criterion (\ref{eqnGIC}) in Section~\ref
{subsecModelSelection}, which can be used to select the true model
among candidate models $\mathrm I$--$\mathrm{\IV}$ and its selection
consistency is provided by Theorem~\ref{thmGIC}. Therefore, the
current paper sheds new light on distinguishing between nonlinear and
nonstationary generalizations of simple linear regression models, and
the results are applied to find appropriate models for short-term and
long-term interest rates.

\section{Technical proofs}\label{secappendix}
We shall in this section provide technical proofs for Theorems \ref
{thmCLTfhatmhat}--\ref{thmGIC}. Because of the time-varying feature
and nonstationarity, the proofs are much more involved than existing
ones for stationary processes. We shall here use techniques of
martingale approximation and $m$-dependent approximation. Let
$\bolds\varepsilon_i = (\bolds\xi_i^\top,\eta_i)^\top$ and
$\bolds{\mathcal{ F}}_i = (\ldots,\bolds\varepsilon
_{i-1},\bolds\varepsilon_i)$ be the corresponding shift process. We
define the projection operator
\[
\mathcal{P}_k \cdot= E(\cdot\mid\bolds{\mathcal{ F}}_k) - E(\cdot\mid
\bolds{\mathcal{ F}}_{k-1}), \qquad k \in
\mathbb Z.
\]
Throughout this section, $C>0$ denotes a constant whose value may vary
from place to place. Let $\alpha_{i,n}(\boldsymbol u,t)$, $i =
1,\ldots
,n$, be a triangular array of deterministic nonnegative weight
functions, $(\boldsymbol u,t) \in\mathbb R^d \times[0,1]$. Lemma~\ref
{lemQbnds} provides a bound for the quantity
\[
Q_\alpha(\boldsymbol u,t) = \sum_{i=1}^n
\bigl\{f_1(\boldsymbol u,i/n \mid\bolds{\mathcal{
F}}_{i-1}) - f(\boldsymbol u,i/n)\bigr\} \alpha_{i,n}(
\boldsymbol u,t),
\]
and is useful for proving Theorems \ref{thmCLTfhatmhat}--\ref{thmunifbndfmu}.

%
\begin{lemma}\label{lemQbnds}
Let $A_n(\boldsymbol u,t) = \max_{1 \leq i \leq n} |\alpha
_{i,n}(\boldsymbol u,t)|$ and define $\bar A_n(\boldsymbol u,t) =
n^{-1} \sum_{i=1}^n |\alpha_{i,n}(\boldsymbol u,t)|$. Then $\|
Q_\alpha
(\boldsymbol u,t)\| \leq\{n A_n(\boldsymbol u,t) \bar A_n(\boldsymbol
u,t)\}^{1/2} \Psi_{0,2}$.
\end{lemma}

\begin{pf}
Since $\mathcal P_k Q_\alpha(\boldsymbol u,t)$, $k \in\mathbb Z$ form
a sequence of martingale differences, we have
\begin{eqnarray*}
\bigl\|Q_\alpha(\boldsymbol u,t)\bigr\|^2 & = & \sum
_{k = -\infty}^n \Biggl\llVert\sum
_{i=1}^n \mathcal P_k \bigl
\{f_1(\boldsymbol u,i/n \mid\bolds{\mathcal{ F}}_{i-1})
\bigr\} \alpha_{i,n}(\boldsymbol u,t)\Biggr\rrVert^2
\\
& \leq& \sum_{k = -\infty}^n \Biggl\{\sum
_{i=1}^n \psi_{i-k-1,2}\bigl |
\alpha_{i,n}(\boldsymbol u,t)\bigr| \Biggr\}^2,
\end{eqnarray*}
and the result follows by observing that $\sum_{i=1}^n \psi_{i-k-1,2}
|\alpha_{i,n}(\boldsymbol u,t)| \leq\break  A_n(\boldsymbol u,t) \Psi_{0,2}$
and $\sum_{i=1}^n \sum_{k \in\mathbb Z} \psi_{i-k-1,2} |\alpha
_{i,n}(\boldsymbol u,t)| \leq n \bar A_n(\boldsymbol u,t) \Psi_{0,2}$.
\end{pf}

%
\begin{lemma}\label{lemCLTThatTtilde}
Assume \textup{(A1)--(A3)} and $\eta_i \in\mathcal L^p$, $p > 2$, $i =
1,\ldots,n$. \textup{(i)} If $b_n \to0$, $h_n \to0$ and $nb_nh_n^d \to
\infty$, then for any $(\boldsymbol u,t) \in\mathbb{R}^d \times
(0,1)$,
\[
\bigl(nb_nh_n^d\bigr)^{1/2} \bigl[
\hat T_n(\boldsymbol u,t) - E\bigl\{\hat T_n(\boldsymbol
u,t)\bigr\}\bigr] \Rightarrow N\bigl[0, \bigl\{m(\boldsymbol u,t)^2 +
\sigma(\boldsymbol u,t)^2\bigr\} f(\boldsymbol u,t)
\lambda_K\bigr],
\]
where $\lambda_K = \lambda_{K_S} \lambda_{K_T}$. \textup{(ii)} If $h_n \to0$
and $nh_n^d \to\infty$, then for any
$\boldsymbol u \in\mathbb{R}^d$,
\[
\bigl(nh_n^d\bigr)^{1/2} \bigl[\tilde
T_n(\boldsymbol u) - E\bigl\{\tilde T_n(\boldsymbol u)
\bigr\}\bigr] \Rightarrow N \biggl[0, \lambda_{K_S} \int
_0^1 \bigl\{m(\boldsymbol u,t)^2 +
\sigma(\boldsymbol u,t)^2\bigr\} f(\boldsymbol u,t) \,dt \biggr].
\]
\end{lemma}

\begin{pf}
Write
\[
\hat T_n(\boldsymbol u,t) - E\bigl\{\hat T_n(\boldsymbol
u,t)\bigr\} = M_n(\boldsymbol u,t) + N_n(\boldsymbol
u,t),
\]
where
\[
M_n(\boldsymbol u,t) = \sum_{i=1}^n
\bigl[y_i K_{S,h_n}(\boldsymbol u - \boldsymbol
x_i) - E\bigl\{y_i K_{S,h_n}(\boldsymbol u -
\boldsymbol x_i) \mid\bolds{\mathcal{ F}}_{i-1}\bigr\}
\bigr]w_{b_n,i}(t)
\]
has summands of martingale differences, and
\[
N_n(\boldsymbol u,t) = \sum_{i=1}^n
\bigl[E\bigl\{y_i K_{S,h_n}(\boldsymbol u - \boldsymbol
x_i) \mid\bolds{\mathcal{ F}}_{i-1}\bigr\} - E\bigl
\{y_i K_{S,h_n}(\boldsymbol u - \boldsymbol x_i)
\bigr\}\bigr]w_{b_n,i}(t)
\]
is the remaining term. Let $\alpha_{i,n}(\boldsymbol u,t) =
m(\boldsymbol u,i/n) w_{b_n,i}(t)$, and by Lemma~\ref{lemQbnds},
\[
\bigl\|N_n(\boldsymbol u,t)\bigr\| \leq\int_{[-1,1]^d}
K_S(\boldsymbol s) \bigl\|Q_\alpha(\boldsymbol u -
h_n \boldsymbol s,t)\bigr\| \,d \boldsymbol s = O\bigl\{(nb_n)^{-1/2}
\bigr\}.
\]
We apply the martingale central limit theorem on $M_n(\boldsymbol u,t)$
to show (i). Since
\begin{eqnarray*}
& & \sum_{i=1}^n \bigl\|\bigl[y_i
K_{S,h_n}(\boldsymbol u - \boldsymbol x_i) - E\bigl
\{y_i K_{S,h_n}(\boldsymbol u - \boldsymbol x_i)
\mid\bolds{\mathcal{ F}}_{i-1}\bigr\}\bigr]w_{b_n,i}(t)
\bigr\|_p^p
\\
& &\qquad\leq \sum_{i=1}^n
2^p \bigl\|y_i K_{S,h_n}(\boldsymbol u - \boldsymbol
x_i)\bigr\|_p^p w_{b_n,i}(t)^p
= O\bigl\{\bigl(nb_nh_n^d
\bigr)^{1-p}\bigr\},
\end{eqnarray*}
the Lindeberg condition is satisfied by observing that $p > 2$. Let
\[
L_n(\boldsymbol s,t) = \sum_{i=1}^n
\bigl\{m(\boldsymbol s,i/n)^2 + \sigma(\boldsymbol
s,i/n)^2\bigr\} \bigl\{f_1(\boldsymbol s,i/n \mid
\bolds{\mathcal{ F}}_{i-1}) - f(\boldsymbol s,i/n)\bigr\}
w_{b_n,i}(t)^2.
\]
Then by (A1) and Lemma~\ref{lemQbnds},
\begin{eqnarray*}
& & h_n^d \sum_{i=1}^n
\bigl[E\bigl\{y_i^2 K_{S,h_n}(\boldsymbol u -
\boldsymbol x_i)^2 \mid\bolds{\mathcal{F}}_{i-1}\bigr\} - E\bigl\{y_i^2
K_{S,h_n}(\boldsymbol u - \boldsymbol x_i)^2\bigr
\}\bigr] w_{b_n,i}(t)^2
\\
&&\qquad =  \int_{[-1,1]^d} K_S(\boldsymbol
s)^2 L_n(\boldsymbol u - h_n \boldsymbol s,t)
\,d \boldsymbol s = O_p\bigl\{(nb_n)^{-3/2}\bigr\}.
\end{eqnarray*}
Also, write $E\{y_i K_{S,h_n}(\boldsymbol u - \boldsymbol x_i) \mid
\bolds{\mathcal{ F}}_{i-1}\} = \int_{[-1,1]^d} m(\boldsymbol u - h_n
\boldsymbol s) K_S(\boldsymbol s)\times\break  f_1(\boldsymbol u - h_n\boldsymbol
s,i/n \mid \bolds{\mathcal{ F}}_{i-1}) \,d\boldsymbol s$. Then we have
\[
\bigl(nb_nh_n^d\bigr) \sum
_{i=1}^n \bigl\|E\bigl\{y_i
K_{S,h_n}(\boldsymbol u - \boldsymbol x_i) \mid\bolds{\mathcal{ F}}_{i-1}\bigr\}\bigr\|^2 w_{b_n,i}(t)^2
= O\bigl(h_n^d\bigr),
\]
and (i) follows by $(nb_nh_n^d) \sum_{i=1}^n E\{y_i^2
K_{S,h_n}(\boldsymbol u -
\boldsymbol x_i)^2\} w_{b_n,i}(t)^2 = \{m(\boldsymbol u,t)^2 +
\sigma(\boldsymbol u,t)^2\} f(\boldsymbol u,t) \lambda_{K_S}
\lambda_{K_T} + o(1)$. Case (ii) can be similarly proved.
\end{pf}

\begin{pf*}{Proofs of Theorems \ref{thmCLTfhatmhat} and \ref{thmCLTftildemuhat}}
Letting $m \equiv1$ and $\sigma\equiv0$ in Lemma~\ref
{lemCLTThatTtilde}, (\ref{eqnCLTfhat}) and
(\ref{eqnCLTftilde}) follow directly. For (\ref{eqnCLTmhat}),
write
\[
\hat T_n(\boldsymbol u,t) - \hat f_n(\boldsymbol u,t)
\frac{E\{\hat
T_n(\boldsymbol u,t)\} }{ E\{\hat f_n(\boldsymbol u,t)\}} = I_n + \mathit{\II}_n,
\]
where
\[
I_n = \bigl[\hat f_n(\boldsymbol u,t) - E\bigl\{\hat
f_n(\boldsymbol u,t)\bigr\}\bigr] \biggl[m(\boldsymbol u,t) -
\frac{E\{\hat T_n(\boldsymbol u,t)\} }{
E\{\hat f_n(\boldsymbol u,t)\}} \biggr] = o_p\bigl\{\bigl(nb_nh_n^d
\bigr)^{-1/2}\bigr\}
\]
and
\[
\mathit{\II}_n = \bigl\{\hat T_n(\boldsymbol u,t) - m(\boldsymbol
u,t) \hat f_n(\boldsymbol u,t)\bigr\} - E\bigl\{\hat T_n(
\boldsymbol u,t) - m(\boldsymbol u,t) \hat f_n(\boldsymbol u,t)\bigr
\}.
\]
Note that
\[
\hat T_n(\boldsymbol u,t) - m(\boldsymbol u,t) \hat f_n(
\boldsymbol u,t) = \sum_{i=1}^n \bigl
\{y_i - m(\boldsymbol u,t)\bigr\} K_{S,h_n}(\boldsymbol u -
\boldsymbol x_i) w_{b_n,i}(t),
\]
by Lemma~\ref{lemCLTThatTtilde}(i),
\[
\bigl(nb_nh_n^d\bigr)^{1/2}
\mathit{\II}_n \Rightarrow N\bigl\{0, \sigma(\boldsymbol u,t)^2 f(
\boldsymbol u,t) \lambda_{K_S} \lambda_{K_T}\bigr\}.
\]
Since $\hat f_n(\boldsymbol u,t) \to f(\boldsymbol u,t)$ in
probability, (\ref{eqnCLTmhat}) follows by Slutsky's theorem.
Case~(\ref{eqnCLTmuhat}) can be similarly proved.
\end{pf*}

\begin{pf*}{Proofs of Theorems \ref{thmunifbnds} and \ref{thmunifbndfmu}} We shall first\vspace*{1pt}
prove Theorem~\ref{thmunifbnds}(i). For this, since $\sup_{t \in
[0,1]} \|\boldsymbol G(t;\bolds{\mathcal{ H}}_0) \|_r < \infty$, we
have $\max_{1 \leq i \leq n} |\boldsymbol x_i| = o_p(n^{1/r'})$ for any
$r' < r$. Hence, $\sup_{t \in[0,1]} \sup_{|\boldsymbol u| > n^{1/r'}}
\hat f_n(\boldsymbol u,t) = 0$ almost surely, and $\sup_{t \in[0,1]}
\sup_{|\boldsymbol u| > n^{1/r'}} E\{\hat f_n(\boldsymbol u,t)\} =
O(n^{-1} h_n^{-d}) = o\{(nb_nh_n^d)^{-1/2}\}$. Therefore, it suffices
to deal with the case in which $|\boldsymbol u| \leq n^{1/r'}$. We
shall here assume that $d = 1$. Cases with higher dimensions can be
similarly proved without extra essential difficulties, but they aew
technically tedious. Let
%
%
\begin{eqnarray}
\label{eqN27204} \hat f^\circ_n(\boldsymbol u,t) & = & \sum
_{i=1}^n E\bigl\{ K_{S,h_n}(
\boldsymbol u - \boldsymbol x_i) w_{b_n,i}(t) \mid\bolds{\mathcal{F}}_{i-1}\bigr\}
\nonumber
\\[-8pt]
\\[-8pt]
\nonumber
& = & \sum_{i=1}^n w_{b_n,i}(t)
\int K_S(\boldsymbol s) f_1(\boldsymbol u -
h_n \boldsymbol s, i/n \mid\bolds{\mathcal{ F}}_{i-1}) \,d
\boldsymbol s.
\end{eqnarray}
Observe that $K_{S,h_n}(\boldsymbol u - \boldsymbol x_i) w_{b_n,i}(t) -
E\{K_{S,h_n}(\boldsymbol u - \boldsymbol x_i) w_{b_n,i}(t) \mid
\bolds{\mathcal{ F}}_{i-1}\}$, $i = 1,\ldots,n$, form a sequence of
bounded martingale differences. By the inequality of \citet
{Freedman1975} and the proof of Theorem~2 in \citet
{WuHuangHuang2010}, we obtain that, for some large constant $\lambda
> 0$,
\[
\mathrm{pr} \Bigl\{\sup_{t \in[0,1]} \sup_{|\boldsymbol u| \leq
n^{1/r'}} \bigl|
\hat f_n(\boldsymbol u,t) - \hat f^\circ_n(
\boldsymbol u,t)\bigr| \geq\lambda(nb_nh_n)^{-1/2} (
\log n)^{1/2} \Bigr\} = o\bigl(n^{-2}\bigr).
\]
Let $\vartheta_i(\boldsymbol u) = f_1(\boldsymbol u, i/n \mid
\bolds{\mathcal{ F}}_{i-1}) - f(\boldsymbol u, i/n)$ and $\Theta
_{l,j}(\boldsymbol u) = \sum_{i=l}^{l+j} \vartheta_i(\boldsymbol u)$.
By (\ref{eqN27204}) and the proof of Lemma~5.3 in \citet
{ZhangWu2012}, it suffices to show that for all $l$,
%
%
\begin{equation}
\label{eqN27226} \mathrm{pr} \Bigl\{\max_{0 \leq j \leq nb_n} \sup
_{|\boldsymbol u| \leq
n^{1/r'}} \bigl|\Theta_{l,j}(\boldsymbol u)\bigr| \geq
\bigl(h_n^{-1} n b_n \log n\bigr)^{1/2}
\Bigr\} = o(b_n).
\end{equation}
Let $\Delta= (nb_nh_n)^{-1/2} (\log n)^{1/4}$ and $\lfloor\boldsymbol
u \rfloor_\Delta= \Delta\lfloor\boldsymbol u / \Delta\rfloor$. By
Theorem~2(ii) in \citet{LiuXiaoWu2013}, under condition (A4),
%
%
\begin{eqnarray}
\label{eqN29948} & & \mathrm{pr} \Bigl\{\max_{0 \leq j \leq nb_n}
\sup
_{|\boldsymbol u|
\leq n^{1/r'}}\bigl |\Theta_{l,j}\bigl(\lfloor\boldsymbol u
\rfloor_\Delta\bigr)\bigr| \geq\bigl(h_n^{-1} n
b_n \log n\bigr)^{1/2} \Bigr\}
\nonumber
\\[-8pt]
\\[-8pt]
\nonumber
&&\qquad =  O \biggl\{\frac{nb_n \Delta^{-1} n^{1/r'} }{(h_n^{-1} n b_n
\log
n)^{q/2}} \biggr\}.
\end{eqnarray}
By (A3), $\max_{0 \leq j \leq nb_n} \sup_{|\boldsymbol u| \leq
n^{1/r'}} |\Theta_{l,j}(\boldsymbol u) - \Theta_{l,j}(\lfloor
\boldsymbol u \rfloor_\Delta) | = O(nb_n \Delta)$, (\ref{eqN27226})
follows. For Theorem~\ref{thmunifbnds}(ii), by Lemma~\ref
{lemunifbnds}, $\sup_{t \in[0,1]} \sup_{\boldsymbol u \in\mathscr X}
|\hat T_n(\boldsymbol u,t) -\break   E\{\hat T_n(\boldsymbol u,t)\}| = O_p\{
(nb_nh_n^d)^{-1/2} (\log n)^{1/2} + (nb_nh_n^d)^{-1}(n^{1/p}\log n)\}$. Since
\begin{eqnarray*}
&&\hat f_n(\boldsymbol u,t) \biggl[\hat m_n(\boldsymbol
u,t) - \frac{E\{\hat
T_n(\boldsymbol u,t)\} }{ E\{\hat f_n(\boldsymbol u,t)\}} \biggr] \\
&&\qquad =\hat T_n(\boldsymbol u,t) - E
\bigl\{\hat T_n(\boldsymbol u,t)\bigr\}
 + E\bigl\{\hat T_n(\boldsymbol u,t)\bigr\} \biggl[1 -
\frac{\hat f_n(\boldsymbol
u,t) }{ E\{\hat f_n(\boldsymbol u,t)\}} \biggr],
\end{eqnarray*}
the result follows. Theorem~\ref{thmunifbndfmu} can be similarly proved.
\end{pf*}

Recall that $\mathscr X \in\mathbb R^d$ is a compact set. Lemma~\ref
{lemunifbnds} provides uniform bounds for
\begin{eqnarray*}
\hat U(\boldsymbol u,t) & = & \sum_{i=1}^n
m(\boldsymbol x_i,i/n) K_{S,h_n}(\boldsymbol u - \boldsymbol
x_i) w_{b_n,i}(t);
\\
\hat V(\boldsymbol u,t) & = & \sum_{i=1}^n
\sigma(\boldsymbol x_i,i/n) \eta_i K_{S,h_n}(
\boldsymbol u - \boldsymbol x_i) w_{b_n,i}(t);
\\
\tilde U(\boldsymbol u) & = & n^{-1} \sum
_{i=1}^n m(\boldsymbol x_i,i/n)
K_{S,h_n}(\boldsymbol u - \boldsymbol x_i);
\\
\tilde V(\boldsymbol u) & = & n^{-1} \sum
_{i=1}^n \sigma(\boldsymbol x_i,i/n)
\eta_i K_{S,h_n}(\boldsymbol u - \boldsymbol
x_i),
\end{eqnarray*}
and is useful in proving Theorems \ref{thmunifbnds} and \ref{thmunifbndfmu}.

%
\begin{lemma}\label{lemunifbnds}
Assume \textup{(A1)}, \textup{(A3)}, \textup{(A4)}, $\eta_i \in\mathcal L^p$ for some $p > 2$, $i
= 1,\ldots,n$, $b_n \to0$ and $h_n \to0$. Let $\chi_n = n^{1/p}\log
n$. \textup{(i)} If $nb_nh_n^d \to\infty$ and $n^{2+d-q} b_n^{d-q} h_n^{d(d+q)}
\to0$, then
%
%
\begin{eqnarray}
\sup_{t \in[0,1]} \sup_{\boldsymbol u
\in\mathscr X} \bigl|\hat U(\boldsymbol
u,t)\bigr| & = & O_p\bigl\{ \bigl(nb_nh_n^d
\bigr)^{-1/2}(\log n)^{1/2}\bigr\}, \label{eqnunifbndUut}
\\
\sup_{t \in[0,1]} \sup_{\boldsymbol u
\in\mathscr X} \bigl|\hat V(\boldsymbol
u,t)\bigr| & = & O_p\bigl\{ \bigl(nb_nh_n^d
\bigr)^{-1/2}(\log n)^{1/2} + \bigl(nb_nh_n^d
\bigr)^{-1}\chi_n\bigr\}. \label
{eqnunifbndVut}
\end{eqnarray}
\textup{(ii)} If $nh_n^d \to\infty$ and $n^{2+d-q} h_n^{d(d+q)} \to0$, then
%
%
\begin{eqnarray}
\sup_{\boldsymbol u \in\mathscr X} \bigl|\tilde U(\boldsymbol u)\bigr| & = & O_p
\bigl\{\bigl(nh_n^d\bigr)^{-1/2}(\log
n)^{1/2}\bigr\}, \label{eqnunifbndUu}
\\
\sup_{\boldsymbol u \in\mathscr X} \bigl|\tilde V(\boldsymbol u)\bigr| & = & O_p
\bigl\{\bigl(nh_n^d\bigr)^{-1/2}(\log
n)^{1/2} + \bigl(nh_n^d\bigr)^{-1}
\chi_n\bigr\}. \label
{eqnunifbndVu}
\end{eqnarray}
\end{lemma}

\begin{pf}
The proof of (\ref{eqnunifbndUut}) is similar to that of Theorem~\ref
{thmunifbnds}(i), and we shall only outline the key differences.
First, the supreme in (\ref{eqnunifbndUut}) is taken over $\boldsymbol
u \in\mathscr X$, a compact set, instead of $\mathbb R^d$. Hence the
truncation argument is no longer needed, and the term $\Delta^{-1}
n^{1/r'}$ in (\ref{eqN29948}) can be replaced by $\Delta^{-1}$.
Second, $E\{m(\boldsymbol x_i,i/n) K_{S,h_n}(\boldsymbol u -
\boldsymbol x_i) \mid\bolds{\mathcal{ F}}_{i-1}\} = \int_{[-1,1]^d}
K_S(\boldsymbol s) f_1^\dagger(\boldsymbol u - h_n \boldsymbol s, i/n
\mid\bolds{\mathcal{ F}}_{i-1}) \,d \boldsymbol s$, where
$f_1^\dagger
(\boldsymbol u, t \mid\bolds{\mathcal{ F}}_{i-1}) = m(\boldsymbol
u, t) f_1(\boldsymbol u,t \mid\bolds{\mathcal{ F}}_{i-1})$. By
(A1), $f_1^\dagger$ satisfies condition (A3), and its predictive
dependence measure is of order (\ref{eqnpsikq}). Hence the proof of
Theorem~\ref{thmunifbnds}(i) applies. Case (\ref{eqnunifbndUu}) can
be similarly handled. For (\ref{eqnunifbndVut}) and (\ref
{eqnunifbndVu}), we shall only provide the proof of (\ref
{eqnunifbndVu}) since (\ref{eqnunifbndVut}) can be similarly derived.
Let $\eta_i^\star= \eta_i \indicator{|\eta_i| \leq n^{1/p}}$ and
$\tilde V^\star(\boldsymbol u)$ be the counterpart of $\tilde
V(\boldsymbol u)$ with $\eta_i$ therein replaced by $\eta_i^\star$, $i
= 1,\ldots,n$. Also, let $\eta_i^\dag= \eta_i^\star- E(\eta
_i^\star
)$, and we can similarly define $\tilde V^\dag(\boldsymbol u)$. Since
$\eta_i \in{\cal L}^p$ are i.i.d., we have $\max_{1 \leq i \leq n}
|\eta_i| = o_p(n^{1/p})$ and $\mathrm{pr}\{\tilde V(\boldsymbol u) =
\tilde V^\star(\boldsymbol u) \mbox{ for all } \boldsymbol u \in
\mathscr X\} \to1$. In addition,
\[
\tilde V^\star(\boldsymbol u) - \tilde V^\dag(\boldsymbol u)
= n^{-1} E\bigl(\eta_i^\star\bigr) \sum
_{i=1}^n \sigma(\boldsymbol x_i,i/n)
K_{S,h_n}(\boldsymbol u - \boldsymbol x_i).
\]
Since $E(\eta_i) = 0$, we have $E(\eta_i^\star) = -E(\eta_i
\indicator
{|\eta_i| > n^{1/p}}) = O(n^{1/p-1})$, and by (\ref{eqnunifbndUu}), it
suffices to show that (\ref{eqnunifbndVu}) holds with $\tilde V^\dag
(\boldsymbol u)$. Let $\bar{\mathscr X} = \{\boldsymbol u \in\mathbb
R^d\dvtx|\boldsymbol u - \boldsymbol v| \leq1 \mbox{ for some }
\boldsymbol v \in\mathscr X\}$, $c_K = \sup_{\boldsymbol v \in
[-1,1]^d} |K_S(\boldsymbol v)|$, $c_1 = \mathrm{var}(\eta_i^\star)$ and
$c_2 = \sup_{t \in[0,1]} \sup_{\boldsymbol u \in\bar{\mathscr X}}
\sigma(\boldsymbol u,t)^2 < \infty$ under (A1). Recall $c_0$ from (A3),
then\break $|\sigma(\boldsymbol x_i, i/n) \eta_i^\dag K_{S,h_n}(\boldsymbol u
- \boldsymbol x_i)| \leq2 c_2^{1/2} c_K n^{1/p} h_n^{-d}$ and
\[
E\bigl\{\sigma(\boldsymbol x_i,i/n)^2 \bigl(
\eta_i^\dag\bigr)^2 K_{S,h_n}(
\boldsymbol u - \boldsymbol x_i)^2 \mid\bolds{\mathcal{ F}}_{i-1}\bigr\} \leq h_n^{-d}c_0c_1c_2
\lambda_{K_S}.
\]
Let $\varpi_n = (nh_n^d)^{-1/2} (\log n)^{1/2} + (nh_n^d)^{-1}(n^{1/p}
\log n)$. Applying the inequality of \citet{Freedman1975} to $\tilde
V^\dag(\boldsymbol u)$, we obtain that, for some large constant
$\lambda> 0$,
\begin{eqnarray*}
& & \mathrm{pr} \bigl\{\bigl|\tilde V^\dag(\boldsymbol u)\bigr| \geq\lambda
\varpi_n\bigr\}
\\
&&\qquad \leq 2 \exp\biggl(-\frac{\lambda^2 \varpi_n^2 }{4c_2^{1/2}c_K
\lambda n^{1/p-1} h_n^{-d} \varpi_n + 2c_0c_1c_2\lambda
_{K_S}n^{-1}h_n^{-d}} \biggr) = O\bigl(n^{-\lambda^{1/2}}
\bigr),
\end{eqnarray*}
and (\ref{eqnunifbndVu}) follows by the discretization argument as in
(\ref{eqN29948}).
\end{pf}

Let $\omega_n = (nb_nh_n^d)^{-1}\log n + b_n^4 + h_n^4$, Lemmas \ref
{lemRSS1}--\ref{lemRSS4} provide asymptotic properties of the
restricted residual sum of squares for models $\mathrm I$--$\mathrm
{\IV
}$, respectively, and are useful in proving Theorem~\ref{thmGIC}. We
shall here only provide the proof of Lemmas~\ref{lemRSS1} and \ref
{lemRSS2}, which relate to nonparametric kernel estimation of
nonlinear regression functions that have been studied in Sections~\ref
{subsecCLTs} and \ref{subsecunifbnds}. Lemmas \ref{lemRSS3} and \ref
{lemRSS4} relate to linear models with time-varying and time-constant
coefficients, and the proof is available in the supplementary material [\citet{ZhangWu2014}].

%
\begin{lemma}\label{lemRSS1}
Assume \textup{(A1)}, \textup{(A3)--(A5)}, $\eta_i \in\mathcal L^p$ for some $p > 2$, $i =
1,\ldots,n$, $b_n \to0$, $h_n \to0$ and $nb_nh_n^d / (\log n)^2 \to
\infty$. If $n^{2+d-q} b_n^{d-q} h_n^{d(d+q)} \to0$ and $n^{1/p-1/2}
b_n^{-1/2} h_n^{-d/2} \to0$, then
\[
n^{-1}\textsc{rss}_n(\mathscr X,\mathscr T,\mathrm I) =
n^{-1}\sum_{i
\in\mathscr I_n} \indicator{\boldsymbol
x_i \in\mathscr X} e_i^2 + O_p
\biggl\{\omega_n +\frac{b_n + h_n }{(nh_n^d)^{1/2}} \biggr\}.
\]
\end{lemma}

\begin{pf}
Note that one can have the decomposition
\[
n^{-1}\textsc{rss}_n(\mathscr X,\mathscr T,\mathrm I) =
n^{-1}\sum_{i
\in\mathscr I_n} \indicator{\boldsymbol
x_i \in\mathscr X} e_i^2 + I_n
- 2\mathit{\II}_n,
\]
where $I_n = n^{-1}\sum_{i \in\mathscr I_n} \{\hat m_n(\boldsymbol
x_i,i/n) - m(\boldsymbol x_i,i/n)\}^2 \indicator{\boldsymbol x_i \in
\mathscr X} = O_p(\omega_n)$ by Theorem~\ref{thmunifbnds}, and
\[
\mathit{\II}_n = n^{-1}\sum_{i \in\mathscr I_n}
\bigl\{\hat m_n(\boldsymbol x_i,i/n) - m(\boldsymbol
x_i,i/n)\bigr\} \indicator{\boldsymbol x_i \in\mathscr
X} e_i.
\]
We shall now deal with the term $\mathit{\II}_n$. By Lemma~\ref{lemunifbnds}(i) and Theorem~\ref{thmunifbnds}, $\sup_{t \in\mathscr T} \sup
_{\boldsymbol u \in\mathscr X} |\{\hat f_n(\boldsymbol u,t) -
f(\boldsymbol u,t)\}\{\hat m_n(\boldsymbol u,t) - m(\boldsymbol u,t)\}|
= O_p(\omega_n)$ and thus
\[
\sup_{t \in\mathscr T} \sup_{\boldsymbol u \in\mathscr X} \biggl
\vert\hat
m_n(\boldsymbol u,t) - m(\boldsymbol u,t) - \frac{\hat T_n(\boldsymbol u,t)
- m(\boldsymbol u,t) \hat f_n(\boldsymbol u,t) }{ f(\boldsymbol
u,t)}\biggr
\vert= O_p(\omega_n).
\]
Let $\Xi_{i,j,n} = \{m(\boldsymbol x_j,j/n) - m(\boldsymbol x_i,i/n)\}
$, and we can then write
\[
\mathit{\II}_n = \mathit{\II}_{n,L} + \mathit{\II}_{n,Q} + O_p(
\omega_n),
\]
where
\[
\mathit{\II}_{n,L} = n^{-1} \sum_{i \in\mathscr I_n}
\frac{\sum_{j=1}^n \Xi_{i,j,n}
K_{S,h_n}(\boldsymbol x_i - \boldsymbol x_j) w_{b_n,j}(i/n) }{
f(\boldsymbol x_i,i/n)} \indicator{\boldsymbol x_i \in\mathscr X}
e_i
\]
and
\[
\mathit{\II}_{n,Q} = n^{-1} \sum_{i \in\mathscr I_n}
\sum_{j=1}^n \frac{K_{S,h_n}(\boldsymbol x_i - \boldsymbol x_j)
w_{b_n,j}(i/n) \indicator
{\boldsymbol x_i \in\mathscr X} }{ f(\boldsymbol x_i,i/n)}
e_i e_j.
\]
Using the orthogonality of martingale differences and Lemma~2 of \citet
{WuHuangHuang2010}, we have $\mathit{\II}_{n,L} = O_p\{(nh_n^d)^{-1/2} (b_n +
h_n)\}$. Also, by splitting the sum in $\mathit{\II}_{n,Q}$ for cases with $i =
j$ and $i \neq j$, one can have\vspace*{2pt} $\mathit{\II}_{n,Q} = O_p\{(nb_n)^{-1} +
n^{-1/2}(nb_nh_n^d)^{-1/2}\}$. Lemma~\ref{lemRSS1} follows by
$(b_nh_n^d)^{-1/2} = \break o\{(b_nh_n^d)^{-1}\}$.
\end{pf}

%
\begin{lemma}\label{lemRSS2}
Assume \textup{(A1)}, \textup{(A3)--(A5)}, $\eta_i \in\mathcal L^p$ for some $p > 2$, $i =
1,\ldots,n$, $h_n \to0$ and $nh_n^d \to\infty$. If $n^{2+d-q}
h_n^{d(d+q)} \to0$ and $n^{1/p-1/2} h_n^{-d/2} \to0$, then \textup{(i)}
\begin{eqnarray*}
n^{-1} \textsc{rss}_n(\mathscr X, \mathscr T, \mathrm{
\II}) & = & \int_{\mathscr X} \int_{\mathscr T} \bigl
\{m(\boldsymbol u,t) - \bar m(\boldsymbol u)\bigr\}^2 f(\boldsymbol
u,t) \,dt \,d \boldsymbol u
\nonumber
\\
& &{} + n^{-1} \sum_{i \in\mathscr I_n}
e_i^2 \indicator{\boldsymbol x_i \in
\mathscr X} + O_p \biggl\{ \biggl(\frac{\log n }{ nh_n^d}
\biggr)^{1/2} + h_n^2 \biggr\}.
\end{eqnarray*}
\textup{(ii)} If in addition model $\mathrm{\II}$ is correctly specified, then
\[
n^{-1} \textsc{rss}_n(\mathscr X, \mathscr T, \mathrm{
\II}) = n^{-1} \sum_{i \in\mathscr I_n}
e_i^2 \indicator{\boldsymbol x_i \in
\mathscr X} + O_p \biggl\{\frac{\log n }{ nh_n^d} + h_n^4
+ \frac{h_n
}{(nh_n^d)^{1/2}} \biggr\}.
\]
\end{lemma}

\begin{pf}
By Theorem~\ref{thmunifbndfmu},
\begin{eqnarray*}
\textsc{rss}_n(\mathscr X, \mathscr T, \mathrm{\II}) & = & \sum
_{i \in
\mathscr I_n} \bigl[\bigl\{y_i - \bar m(
\boldsymbol x_i)\bigr\} - \bigl\{\hat\mu_n(\boldsymbol
x_i) - \bar m(\boldsymbol x_i)\bigr\}
\bigr]^2 \indicator{\boldsymbol x_i \in\mathscr X}
\\
& = & I_n + O_p\bigl[n\bigl\{\bigl(nh_n^d
\bigr)^{-1/2} (\log n)^{1/2} + h_n^2\bigr
\}\bigr],
\end{eqnarray*}
where by Lemma~2 in \citet{WuHuangHuang2010},
\begin{eqnarray*}
I_n & = & \sum_{i \in\mathscr I_n} \bigl[\bigl
\{y_i - m(\boldsymbol x_i,i/n)\bigr\} + \bigl\{m(
\boldsymbol x_i,i/n) - \bar m(\boldsymbol x_i)\bigr\}
\bigr]^2 \indicator{\boldsymbol x_i \in\mathscr X}
\\
& = & \sum_{i \in\mathscr I_n} \bigl\{m(\boldsymbol
x_i,i/n) - \bar m(\boldsymbol x_i)\bigr\}^2
\indicator{\boldsymbol x_i \in\mathscr X} + \sum
_{i \in\mathscr I_n} e_i^2 \indicator{\boldsymbol
x_i \in\mathscr X} + O_p\bigl(n^{1/2}\bigr).
\end{eqnarray*}
Since $\mathscr X \in\mathbb{R}^d$ is a compact set, by the proof
of Lemma~\ref{lemCLTThatTtilde}, we have
\begin{eqnarray*}
& & \sum_{i \in\mathscr I_n} \bigl\{m(\boldsymbol
x_i,i/n) - \bar m(\boldsymbol x_i)\bigr\}^2
\indicator{\boldsymbol x_i \in\mathscr X}
\\
& &\qquad = \sum_{i \in\mathscr I_n} E\bigl[\bigl\{m(\boldsymbol
x_i,i/n) - \bar m(\boldsymbol x_i)\bigr\}^2
\indicator{\boldsymbol x_i \in\mathscr X}\bigr] + O_p
\bigl(n^{1/2}\bigr)
\\
&&\qquad =  n \int_{\mathscr X} \int_{\mathscr T} \bigl\{m(
\boldsymbol u,t) - \bar m(\boldsymbol u)\bigr\}^2 f(\boldsymbol u,t)
\,dt \,d \boldsymbol u + O\bigl(1 + n^{1/2}\bigr),
\end{eqnarray*}
and (i) follows. Case (ii) follows by a similar argument as in
Lemma~\ref{lemRSS1}.
\end{pf}

%
\begin{lemma}\label{lemRSS3}
Assume \textup{(A1)--(A3)}, \textup{(A6)}, \textup{(P2)} and $\eta_i \in\mathcal L^p$ for some $p
> 2$, $i = 1,\ldots,n$. If $b_n \to0$ and $nb_n \to\infty$, then \textup{(i)}
\begin{eqnarray*}
n^{-1} \textsc{rss}_n(\mathscr X, \mathscr T, \mathrm{
\III}) & = & \int_{\mathscr X} \int_{\mathscr T} \bigl
\{m(\boldsymbol u,t) - \boldsymbol u^\top\bolds
\beta_n(t)\bigr\}^2 f(\boldsymbol u,t) \,dt \,d \boldsymbol u
\\
& &{} + n^{-1} \sum_{i \in\mathscr I_n}
e_i^2 \indicator{\boldsymbol x_i \in
\mathscr X} + O_p(\phi_n).
\end{eqnarray*}
\textup{(ii)} If in addition model $\mathrm{\III}$ is correctly specified, then
\[
n^{-1} \textsc{rss}_n(\mathscr X, \mathscr T, \mathrm{
\III}) = n^{-1} \sum_{i \in\mathscr I_n}
e_i^2 \indicator{\boldsymbol x_i \in
\mathscr X} + O_p \biggl(\phi_n^2 +
\frac{b_n^2 }{
n^{1/2}} \biggr).
\]
\end{lemma}

%
\begin{lemma}\label{lemRSS4}
Assume \textup{(A1)--(A3)}, \textup{(A6)}, \textup{(P1)} and $\eta_i \in\mathcal L^p$ for some $p
> 2$, $i = 1,\ldots,n$. Then \textup{(i)}
\begin{eqnarray*}
n^{-1} \textsc{rss}_n(\mathscr X, \mathscr T, \mathrm{
\IV}) & = & \int_{\mathscr X} \int_{\mathscr T} \bigl
\{m(\boldsymbol u,t) - \boldsymbol u^\top\bolds\theta_n
\bigr\}^2 f(\boldsymbol u,t) \,dt \,d \boldsymbol u
\nonumber
\\
& & {} + n^{-1} \sum_{i \in\mathscr I_n}
e_i^2 \indicator{\boldsymbol x_i \in
\mathscr X} + O_p\bigl(n^{-1/2}\bigr).
\end{eqnarray*}
\textup{(ii)} If in addition model $\mathrm{\IV}$ is correctly specified, then
\[
n^{-1} \textsc{rss}_n(\mathscr X, \mathscr T, \mathrm{
\IV}) = n^{-1} \sum_{i \in\mathscr I_n}
e_i^2 \indicator{\boldsymbol x_i \in
\mathscr X} + O_p\bigl(n^{-1}\bigr).
\]
\end{lemma}

\begin{pf*}{Proof of Theorem~\ref{thmGIC}} For model $\mathrm I$, the AMSE optimal
bandwidths satisfy $b_n(\mathrm I) \asymp n^{-1/(d+5)}$ and
$h_n(\mathrm I) \asymp n^{-1/(d+5)}$. By Lemma~\ref{lemRSS1}, we
have
\[
\log\bigl\{\textsc{rss}_n(\mathscr X, \mathscr T, \mathrm I)/n\bigr
\} = \log\biggl(n^{-1} \sum_{i \in\mathscr I_n}
e_i^2 \indicator{\boldsymbol x_i \in
\mathscr X} \biggr) + O_p\bigl\{n^{-7/(2d+10)}\bigr\}.
\]
Under the stated conditions on the tuning parameter, we have
$n^{-7/(2d+10)} = o\{\tau_n (b_n h_n^d)^{-1}\}$, and thus the
estimation error is dominated by $\tau_n \textsc{df}(\mathrm I)$ which
goes to zero as $n \to\infty$. By Lemmas \ref{lemRSS2}--\ref
{lemRSS4}, similar results can be derived for models $\mathrm{\II
}$--$\mathrm{\IV}$. Note that
\[
\tau_n \max\bigl\{\textsc{df}(\mathrm I), \textsc{df}(\mathrm{\II}),
\textsc{df}(\mathrm{\III}), \textsc{df}(\mathrm{\IV})\bigr\} = o(1),
\]
which will be dominated by any model misspecification. The result
follows by $\textsc{df}(\mathrm{\IV}) < \min\{\textsc{df}(\mathrm
{\II
}), \textsc{df}(\mathrm{\III})\} \leq\max\{\textsc{df}(\mathrm
{\II}),
\textsc{df}(\mathrm{\III})\} < \textsc{df}(\mathrm I)$.
\end{pf*}

\section*{Acknowledgments}
We are grateful to the Editor, an Associate Editor, and two anonymous
referees for their helpful comments and suggestions.


\begin{supplement}[id=suppA]
\stitle{Additional technical proofs}
\slink[doi]{10.1214/14-AOS1299SUPP} 
\sdatatype{.pdf}
\sfilename{aos1299\_supp.pdf}
\sdescription{This supplement contains technical proofs of Lemmas \ref
{lemRSS3} and \ref{lemRSS4}, Proposition~\ref{propS10828a} and
Theorems \ref{thmextensioncorrelation} and \ref{thmextensionautoregression}.}
\end{supplement}

%
%

%

\printaddresses
\end{document}